\documentclass[12pt]{article}

\usepackage{amsthm}
\usepackage{amsmath,amsrefs,amssymb,latexsym,amsfonts,mathrsfs}
\usepackage{graphicx}

\usepackage{color}

\usepackage{fancyhdr}
\usepackage[top=30truemm, bottom=30truemm, left=25truemm, right=25truemm]{geometry}

\usepackage{hyperref} 
\hypersetup{
    colorlinks=true, 
    linktoc= page,   
citecolor=blue
}

\newcommand{\ep}{\varepsilon}

\newcommand{\SCR}[1]{{\mathscr #1}}

\newcommand{\CAL}[1]{{\cal #1}
}

\newcommand{\MAT}[1]{\begin{pmatrix}#1\end{pmatrix}}

\newcommand{\J}[1]{\left\langle #1 \right\rangle}

\theoremstyle{definition}
\newtheorem{Thm}{{\bf Theorem}}[section]

\newtheorem{Lem}[Thm]{{\bf Lemma}}
\newtheorem{Prop}[Thm]{{\bf Proposition}}
\newtheorem{Cor}[Thm]{{\bf Corollary}}
\newtheorem{Ass}[Thm]{{\bf Assumption}}

\newtheorem{Def}[Thm]{{\bf Definition}}
\newtheorem{Rem}[Thm]{{\bf Remark}}

\newcounter{Exami}

\newcommand{\Proof}[2][Proof]{
\begin{proof}[{\bf #1}]
#2
\end{proof}
}

\numberwithin{equation}{section}

\usepackage[italicdiff]{physics}

\begin{document}

\begin{flushleft}
{\bf \Large Global well-posedness and scattering in weighted space for nonlinear Schr\"{o}dinger equations below the Strauss exponent without gauge-invariance
 } \\ \vspace{0.3cm}
 \end{flushleft}
\begin{flushleft}
{\large Masaki KAWAMOTO}\\
{Research Institute for Interdisciplinary Science, Okayama University, Okayama City Okayama, 700-8530, Japan }\\
Email: {kawamoto.masaki@okayama-u.ac.jp}
\end{flushleft}
\begin{flushleft}
{\large Satoshi MASAKI}\\
{Department of mathematics, Hokkaido University, Sapporo Hokkaido, 060-0810, Japan }\\
Email: {masaki@math.sci.hokudai.ac.jp}
\end{flushleft}
\begin{flushleft}
{\large Hayato MIYAZAKI}\\
{
Teacher Training Courses, Faculty of Education, Kagawa University, Takamatsu, Kagawa 760-8522, Japan
\\
Email: miyazaki.hayato@kagawa-u.ac.jp
}
\end{flushleft}
\begin{abstract}
In this paper, we consider the nonlinear Schr\"{o}dinger equation (NLS) with a general homogeneous nonlinearity in dimensions up to three. We assume that the degree (i.e., power) of the nonlinearity is such that the equation is mass-subcritical and short-range. We establish global well-posedness (GWP) and scattering for small data in the standard weighted space for a class of homogeneous nonlinearities, including non-gauge-invariant ones. Additionally, we include the case where the degree is less than or equal to the Strauss exponent. When the nonlinearity is not gauge-invariant, the standard Duhamel formulation fails to work effectively in the weighted Sobolev space; for instance, the Duhamel term may not be well-defined as a Bochner integral. To address this issue, we introduce an alternative formulation that allows us to establish GWP and scattering, even in the presence of poor time continuity of the Duhamel term.
\end{abstract}

\begin{flushleft}

{\em Keywords and phrases:} nonlinear Schr\"{o}dinger equation; global in time well-posedness; nonlinear scattering;  Strauss exponent. \\
{\em Mathematical Subject Classification:} Primary: 35Q55, Secondly: 35B40, 35P25.

\end{flushleft}

\section{Introduction}

In this paper, we consider the Cauchy problem of
the nonlinear Schr\"{o}dinger equation with a homogeneous nonlinearity:
\begin{align} \label{1}
\begin{cases}
& {\displaystyle  i \partial _t u (t,x) = \frac{p^2  u (t,x)}2+F(u(t,x))}, \\
& u(0,x) = u_0(x),
\end{cases}
\end{align}
where $t \in \mathbb{R}$, $x \in \mathbb{R}^d $, $d=1,2, 3$, $p = - i \nabla $, $p^2 = |p|^2 = -\Delta$ and $F$
is homogeneous of degree $\alpha$, i.e.,
\[
	F( ru) = r^{\alpha} F(u)
\]
holds for any $u \in \mathbb{C}$ and $r>0$.
A typical example of a homogeneous nonlinearity is the gauge invariant power type nonlinearity $F(u)=\lambda |u|^{\alpha-1} u$.
The problem arising herein is whether the global-in-time wellposedness (GWP) holds or not when the exponent $\alpha$ is smaller than or equal to the Strauss exponent
\begin{align*}
\alpha_S = 1+ \frac{2-d + \sqrt{d^2 + 12 d + 4} }{2d}.
\end{align*}
Moreover, we consider the scattering problem for small data.

As is well known, if $\alpha_S < \alpha < \alpha^{\ast}$, where $\alpha^{\ast} = \infty$ if $d \leq 2$, otherwise, $\alpha^{\ast}= 1 + 4/(d-2)$,
then the GWP and scattering of \eqref{1} follows for small data, in a suitable topology, for a wide class of homogeneous nonlinearities
by means of the standard perturbation argument with exotic Strichartz estimates, see e.g. Kato \cite{Ka2}.
However, if $\alpha$ is smaller than or equal to the Strauss exponent, the GWP and scattering for non-gauge-invariant nonlinearities is not clear even for small data. Indeed, the finite-time blowup is studied for the modulus nonlinearity $F(u)=|u|^\alpha$ (see, e.e., \cite{II,IW,MM2,MS,FO}). As for the equations with the gauge-invariant nonlinearity,
the analysis in weighted spaces, radial negative order Sobolev spaces, or Fourier-Lebesgue spaces
has been employed to obtain the local and global well-posedness and scattering (See \cite{CW1, GOV, NO, NP09, HT, HNT1, HNT2, Ma1, Ma2, Ma3, Hi, GW} and references therein).

The aim of this paper is to conduct an analysis of NLS equation with a broad class of homogeneous nonlinearity, which includes non-gauge-invariant nonlinearities, with exponents smaller than or equal to the Strauss exponent in the framework of the standard weighted $L^2$-space.
We discuss the nontriviality in more detail later in Section \ref{ss:difficulty}.

\subsection{Main result}

We parameterize the power  $\alpha$ of the nonlinearity with $\alpha_0 \in [1, 2]$ using the formula:
\[
	\alpha= \alpha(\alpha_0) := 1+ \frac{2\alpha_0}d
\]
or, equivalently, $\alpha_0=\alpha_0(\alpha) = \frac{d}2 (\alpha-1)$.
Note that $\alpha(1)=1+2/d$ represents the threshold between the short-range case and the long-range case, while $\alpha(2)=1+4/d$ corresponds to the mass-critical exponent. Essentially, we use $\alpha$ to indicate the power of the nonlinearity, and $\alpha_0$ is employed to describe various exponents in norms and other contexts.

Throughout the paper,
we assume
\begin{equation}\label{IC}
\alpha_0\in \begin{cases}
(\frac54,2) & d=1, \\
(1,2) & d=2 , \\
\{ \frac32 \} & d=3,
\end{cases}
\end{equation}
or equivalently
\begin{equation}\label{IC2}
\alpha\in \begin{cases}
(\frac72,5) & d=1, \\
(2,3) & d=2, \\
\{ 2 \} & d=3.
\end{cases}
\end{equation}
We note that the above range includes
the case $\alpha \le \alpha_S$.
Indeed, $\alpha_S(1)=\frac{3 + \sqrt{17}}{2} > 3.561 > 7/2$,
$\alpha_S(2)= 1 + \sqrt{2}>2$, and $\alpha_S(3)=2$.
We emphasize that our range covers whole intercritical range
in two dimensions.

To investigate the structure of the homogeneous nonlinearities,
let us recall an expansion of the nonlinearity introduced in \cite{MM}. A homogeneous nonlinearity $F$ of degree $\alpha $ is written as
\[
	F(u) = \sum_{n \in \mathbb{Z}} \lambda _n  F_n(u)
\]
under a suitable regularity assumption on the function $\theta \mapsto F(e^{i\theta})$,
where $$F_n(u) :=|u|^{\alpha -n} u^n $$ and
\[
	\lambda_n := \frac1{2\pi} \int_0^{2\pi} e^{-in\theta} F(e^{i\theta}) d\theta.
\]
By the expansion, a homogeneous nonlinearity $F$ can be identified with the coefficients $\{ \lambda_n \}_{n\in \mathbb{Z}} \subset \mathbb{C}$.
This enables us to investigate the structure of the nonlinearity in the language of the coefficients.
We remark that the gauge-invariant nonlinearity $F(u)=\lambda |u|^{\alpha-1}u$ corresponds to $\lambda_n = \lambda \delta_{n1}$ and
the modulus-type nonlinearity $F(u)=\lambda |u|^{\alpha}$ to $\lambda_n = \lambda \delta_{n0}$, where $\delta_{ij}$ is the Kronecker delta.

\begin{Ass}\label{A1}
We suppose that $F$ satisfies \eqref{IC2} and
the following two conditions:
\begin{itemize}
\item[(A1)]
 $  \sum_{n\in {\mathbb Z} } \J{n}^{\frac72} |\lambda_n| < \infty $;
\item[(A2)] $
	\lambda_n = 0$ for $ n \le 0$.
\end{itemize}
\end{Ass}
The assumption (A1) is the weighted summability of the coefficients. It is related with the regularity property of the nonlinearity $F$.
The assumption (A2) is the absence of nonpositive coefficients.
It is known that the NLS with
the modulus-type nonlinearity $F=F_0(u) = |u|^{\alpha}$ admits $L^2$-solutions which blow up in finite time (see, \cite{II,IW,MM2,MS, FO}).
The condition (A2) assures the absence of the $F_0$-component.
Note that, however, these blowing up solutions are out of our framework and hence it is not clear whether (A2) or the condition $\lambda_0=0$ is necessary for our result.

\begin{Rem}
We remark that any finite  linear combination of $\{ F_n \}_{n\ge1}$
satisfies (A1) and (A2).
A
non-trivial example of the nonlinearity satisfying (A1) and (A2), which consists of an infinitely sum,
is
\begin{align*}
F(u) = \frac{|u|^{\alpha} u }{|u| -au} =  \sum_{n=1}^{\infty} a^{n-1} F_n(u)  , \quad -1 < a < 1 .
\end{align*}
We provide one more example, which is close to the borderline case. For a positive integer $k$ and $z \in \mathbb{C} $ with $0<|z|<1$, we define
\[
	h_k(z) := \left(\left( \frac{z-1}{z} \right)^k \log \frac{1}{1-z} +\sum_{m=1}^k \frac1m  \left( \frac{z-1}{z} \right)^{k-m} \right)
	=
	 \sum_{n=1}^\infty \frac{k!}{n(n+1)\cdots (n+k)} z^n.
\]
Note that $h_k$ can be continuously extended to $\{ |z| \le 1\}$.
In particular, $h_k(1)=k^{-1}$ and $h_k(0)=0$.
Then, the nonlinearity
\[
	F(u):= |u|^\alpha h_k\left(\frac{u}{|u|}\right) = |u|^\alpha \left(\left( \frac{u-|u|}{u} \right)^k \log \frac{|u|}{|u|-u} +\sum_{m=1}^k \frac1m  \left( \frac{u-|u|}{u} \right)^{k-m} \right), \quad u\neq0,
\]
with $F(0)=0$,
satisfies (A2) for all $k\ge1$, but (A1) holds only for $k\ge 4$ since we have the expansion
\[
	F(u) =  \sum_{n=1}^\infty \frac{k!}{n(n+1)\cdots (n+k)} F_n(u).
\]
The series expansion of $h_k$ is justified, for instance, by the recurrence relation
$h_k(z)-\frac{h_k(z)}z=-\frac1{k+1} + h_{k+1}(z)$ with the convention $h_0(z)=\log \frac1{1-z}$.
\end{Rem}

We introduce function spaces.
Set a space $H^{ a,b } := \{ u \in \SCR{S}'(\mathbb{R}^d) \, ; \,  \| u \|_{H^{a,b}}    < \infty    \} $ with $\| u \|_{H^{a,b}} := \| \J{p}^{a} u \|_{L^2} + \| \J{x} ^b u \|_{L^2}$ for $a,b \in \mathbb{R}$.
We abbreviate $H^{a}=H^{a,0}$.
For $r\in (1,\infty)$ and $s\in \mathbb{R}$,
$H^{s}_r$ is the Sobolev space given by the norm $\|u\|_{H^s_r}:=\|\ev{p}^su\|_{L^r}$.
$U(t)$ stands for the standard  Schr\"odinger group: $U(t) = e^{-it\frac{p^2}2}$. Let $x^2 = |x|^2 $.
We introduce the operator
\[
J(t) := U(t) x U(-t) = - e^{ix^2/2t} t p e^{-ix^2/2t}
= x - tp ,
\]
where the second formula makes sense for $t\neq0$.
Note that the relation
$U(t_1)J(t_2) = J(t_1+t_2)U(t_1)$ holds for any $t_1,t_2 \in \mathbb{R}$.
For a positive integer $k$ and a parameter $t\in \mathbb{R}$,
we define the function space $X^k(t)$ as follows:
\[
	X^k(t) := \{ u\in L^2 ; U(-t) u \in H^{0,k} \} = \{ u\in L^2 ; J(t)^k u \in L^2 \} .
\]
For any $t_1,t_2\in\mathbb{R}$ with $t_1\neq t_2$, it holds that $X^1(t_1)\neq X^1(t_2)$ and $X^1(t_1) \cap X^1(t_2)=H^{1,1}$. Additionally, $X^1(t) \cap H^1 = H^{1,1}$ holds for any $t\in\mathbb{R}$.
It follows that $U(t_1)X(t_2) = X(t_1+t_2) U(t_1)$.
Further, noting that $\overline{U(t) f}= U(-t)\overline{f}$, one sees that, for any $t\in \mathbb{R}$,
 $f \in X^1(t)$ if and only if $\overline{f} \in X^1(-t)$.

Our definition of a solution is as follows:
\begin{Def}\label{D:sol}
Let $I \ni 0$ be an interval.
Let $t_0 \in \mathbb{R}$. Given $u_0 \in X^1(t_0)$,
we say a function $u(t,x)$ is a solution to \eqref{1} on $I$ if $u \in L^{q_0}_{t, \mathrm{loc}}(I;L^{r_0}_x(\mathbb{R}^d))$ satisfies
$U(-t-t_0)u(t) \in C(I;H^{0,1})$ and the identity
\begin{equation}\label{E:I1}
	U(-t-t_0)u(t) = U(-t_0)u_0 - i \int_0^t U(-t_0-s) F(u(s)) ds
\end{equation}
holds for $t\in I$, where
\begin{equation}\label{E:defq0r0}
	(q_0,r_0) = \begin{cases}
	(12,3) & d=1,\\
	(\frac{2}{\alpha_0-1},\frac2{2-\alpha_0}) & d=2, \\
	(4,3) & d=3
	\end{cases}
\end{equation}
is an admissible pair.
\end{Def}

We remark that the notion of a solution depends on the parameter $t_0$, which is specified by the data $u_0$ via the validity of $u_0 \in X^1(t_0)$.
Another remark is that it is not clear whether or not
the Duhamel term (the second integration term) in the right hand side of \eqref{E:I1} makes sense as a Bochner integral in $H^{0,1}$ from the properties $u \in L^{q_0}_{t, \mathrm{loc}}(I;L^{r_0}_x(\mathbb{R}^d))$ and
$U(-t-t_0)u(t) \in C(I;H^{0,1})$.
Nevertheless, the integral makes sense, for instance, as a Bochner integral in $L^{2}$
if $t_0 \neq 0$ and $-t_0 \not\in I$.
Similarly, if $u_0 \in H^1 \cap H^{0,1}$ then the integral is defined as a Bochner integral in $H^1 \cap H^{0,1}$.
Our definition of a solution includes that a function
$I \ni t\mapsto \int_0^t U(-t_0-s) F(u(s)) ds \in H^{0,1}
$ is continuous, thanks to \eqref{E:I1}.
Note that  the integral
makes sense as a Bochner integral in $H^{0,1}$
if and only if this function is \emph{absolutely continuous} in the $H^{0,1}$-topology.

Our main result is the global existence and scattering in one time direction for small data.
\begin{Thm}[Small data scattering in one-time direction] \label{T1}
Let $d=1$, $2$, $3$.
Suppose \eqref{IC2}, (A1), and (A2).
There exists $\ep_0=\ep_0(d,\alpha, F) >0$ such that if
$u_0$ satisfies
 $u_0 \in X^1(t_0)$ and
\begin{align}\label{E:smallcond}
|t_0|^{\frac1{\alpha-1} - \frac{d}4}
(\| u_0 \|_{L^2 }+|t_0|^{-\frac{1}2}\|J(t_0)u_0\|_{L^2}) \le \ep_0
\end{align}
for some  $t_0\neq 0$.
Then, the following two assertions hold.
\begin{itemize}
\item If $t_0>0$ then
there exists a unique solution $u(t)$
on $[0,\infty)$.
Furthermore, there exists $u_+ \in H^{ 0, 1}$ such that
\begin{align*}
\lim_{t \to \infty} \left\|
U(-t-t_0) u(t, \cdot ) - u_+
\right\|_{H^{0, 1}}  = 0.
\end{align*}
\item If $t_0<0$ then
there exists a unique solution $u(t)$
on $(-\infty,0]$.
Furthermore, there exists $u_- \in H^{ 0, 1}$ such that
\begin{align*}
\lim_{t \to -\infty} \left\|
U(-t-t_0) u(t, \cdot ) - u_-
\right\|_{H^{0, 1}}  = 0.
\end{align*}
\end{itemize}
\end{Thm}
The conclusion is the standard 
global existence and scattering in both the positive and negative time direction for small data.
We remark that, however, when $F$ is not a gauge-invariant nonlinearity, it is even highly nontrivial that the right hand side of \eqref{E:I1} belongs to $C(I;H^{0,1})$.
We discuss this issue in more detail in the forthcoming subsection.
\begin{Rem}\label{R:red}
\begin{enumerate}
\item
By the time reversal symmetry, the case $t_0<0$ follows from the case $t_0>0$.
Indeed, if $u_0 \in X^1(t_0)$ then $\overline{u_0} \in X^1(-t_0)$, as mentioned above.
One verifies that
$u(t,x)$ is a solution to \eqref{1} on an interval $I$ such that $u(0)=u_0 \in X^1(t_0)$ if and only if $v(t,x):=\overline{u(-t,x)}$ is also a solution to \eqref{1} on $-I$ such that $v(0)=\overline{u_0} \in X^1(-t_0)$.
If the theorem is true for positive parameter then it
shows the existence $v(t,x)$ on $[0,\infty)$ and the scattering as $t\to \infty$ under the corresponding smallness assumption.
Then, we see that $u(t,x)$ exists on $(-\infty,0]$ and scatters as $t\to -\infty$.
\item By the scaling symmetry,
the theorem is further reduced to the case $t_0=1$.
Indeed, for any $t_0>0$,
$u$ is a solution on $[0,\infty)$ for the data
$u_0\in X^1(t_0)$ if and only if $v(t,x) := t_0^{\frac1{\alpha-1}} u(t_0 t, \sqrt{t_0}x)$ is a solution on $[0,\infty)$ for the data
$v_0(x):=t_0^{\frac1{\alpha-1}}u_0(\sqrt{t_0}x) \in X^1(1)$. It can be verified, for instance, from
the identity
$J(t+1)v(t,x) = t_0^{\frac1{\alpha-1}-\frac12} (J(t_0t+t_0) u(t_0t,\cdot))(\sqrt{t_0}x)$.
If the theorem is true for $t_0=1$
then we obtain a solution $v(t)$ on an interval $[0,\infty)$
in the sense of Definition \ref{D:sol}.
Then, one sees that the function $u(t,x)$ is a desired solution in the sense of Definition \ref{D:sol}.
\item
However, we do not have a result for data in $X^1(0) = H^{0,1}$.
This shows a contrast in the case of the gauge-invariant nonlinearity. We discuss this issue later.
Here we note that
the key difference between $X^1(t_0)$ ($t_0\neq0$) and $X^1(0)$ is
that the validity of the inclusion
\[
	X^1(t_0) \subset
	 \begin{cases}
	L^2 \cap L^\infty & d=1,\\
	\bigcap_{r\in [2,\infty)} L^r & d=2,\\
	L^2 \cap L^6 & d=3
\end{cases}
\]
for $t_0 \neq 0$,
while $X^1(0) \not\subset L^r$ for any $r>2$.
\end{enumerate}
\end{Rem}

As mentioned above, we have $X^1(1) \cap X^1(-1) = H^{1,1}$.
The local well-posedness for arbitrary $H^{1,1}$ data follows by a standard argument (see \cite{Ca}).
We have the following corollary of  Theorem \ref{T1}, which is scattering for both time direction
for data small in $H^{1,1}$.
\begin{Cor}[Small data scattering in $H^{1,1}$]\label{C:both}
Let $d=1$, $2$, $3$.
Suppose \eqref{IC2}, (A1), and (A2).
There exists $\widetilde{\ep}_0=\widetilde{\ep}_0(d,\alpha, F) >0$ such that if
 $u_0 \in H^{1,1}$ satisfies $\| u_0 \|_{H^{1,1} } \le \widetilde{\ep}_0 $ then there exists a unique solution $u(t)$ of \eqref{1} on $\mathbb{R}$.
Furthermore, there exists $u_\pm \in H^{ 1, 1}$ such that
\begin{align*}
\lim_{t \to \pm \infty} \left\|
U(-t) u(t, \cdot ) - u_\pm
\right\|_{H^{1, 1}}  = 0.
\end{align*}
\end{Cor}

\begin{Rem}
Hayashi and Naumkin \cite{HN} studied the same problem in one dimension for nonlinearities expressed as a finite combination of $\lambda_j |u|^{\alpha_j - j} u^j$ for $j \neq 0$ and $\alpha_j > 3$, and established a similar result. In the two-dimensional case, Germain, Masmoudi, and Shatah \cite{GMS2} obtained global well-posedness (GWP) when a specific linear function is multiplied by a nonlinear term. However, to the best of the authors' knowledge, no results exist for equations with a large class of homogeneous nonlinearities in the two-dimensional case.
The three-dimensional quadratic case has been studied by Kawahara \cite{Kawa05, Kawa07} and Germain, Masmoudi, and Shatah \cite{GMS1}. In \cite{GMS1}, the scattering is obtained in $L^2$ for small data in $X^2(t_0)$, $t_0\neq0$, under the condition $\lambda_{n} = 0$ for $n \notin \{-2, 2\}$, i.e., $F(u) = \lambda_2 u^2 + \lambda_{-2} \overline{u}^2$. In contrast, Theorem \ref{T1} and Corollary \ref{C:both} improve upon these results in the following ways: our results allow us to consider a wider class of initial data $X^1(t_0)$, scattering is obtained in the same space as the initial data, and the nonlinearity is permitted to be an infinite sum of terms of the form $\lambda_j |u|^{2 - j} u^j$ for $j \geq 1$.
On the other hand, we require an additional restriction, $\lambda_{-2} = 0$, meaning we cannot handle the nonlinearity $\overline{u}^2$.
\end{Rem}

\subsection{Difficulty with a standard Duhamel formula}\label{ss:difficulty}

The key point of our results lies in the solvability in weighted spaces despite the presence of non-gauge-invariant nonlinear part, i.e., $F_n$ with $n\neq1$, in the expansion of the nonlinearity $F$. We elaborate on explaining why this is nontrivial.

Firstly, let us review the properties of weighted spaces.
By the properties of the weighted spaces mentioned above, given
 solution $u(t)$ of the linear Schr\"odinger equation with
an initial value $u_0 \in X^1(t_0)\setminus H^{1,1}$, $u(t) \in X^1(s)$ holds for $t\in \mathbb{R}$ if and only if $s=t+t_0$. In other words, the function space to which the solution belongs changes with time, and furthermore, these spaces become different over time.

It is well-known that this property holds true also under the evolution of the NLS with the gauge-invariant nonlinearity.
 It stems, for instance, from the fact that $J(t)$ acts like a derivative for the gauge-invariant nonlinearity.
Indeed,
for the gauge-invariant nonlinearity $F_1(u)$, we have
\begin{align*}
	J(t+1) F_1(u)
	={}& - e^{ix^2/2(t+1)} J (t+1) p e^{-ix^2/2(t+1)} F_1(u) \\
	={}& - e^{ix^2/2(t+1)} J (t+1) p  F_1(e^{-ix^2/2(t+1)} u) \\
	={}& \frac{\partial F_1}{\partial z}(u) (J(t+1)u) - \frac{\partial F_1}{\partial \overline{z}}(u) (\overline{J(t+1)u}).
\end{align*}
If the nonlinearity is $F=F_1$ then, using the above identity and the Duhamel formula, we obtain
\[
	J(t+1)u(t) = U(t)J(1) u_0 -i \int_0^t U(t-s) \left( \frac{\partial F_1}{\partial z}(u) (J(\cdot+1)u) - \frac{\partial F_1}{\partial \overline{z}}(u) (\overline{J(\cdot+1)u}) \right)(s) ds.
\]
Using this formula, a solution $u(t)$ with $u(0)=u_0 \in X^1(1)$
in the sense of Definition \ref{D:sol}
can be obtained by a standard fixed point argument
 (see \cite{Ma1,Ma2,NO, HNT1, HNT2} and references therein).
Note that the solution satisfies
 $J(t+1)u(t) \in L^2$ for each $t$.

This discussion also highlights the nontriviality of obtaining local solutions in weighted spaces in our case. In fact, when considering a nonlinear term $F_n$ ($n\neq 0,1$), the construction of a solution for $u(0)=u_0 \in X^1(1)$ seems difficult, at least with such a direct argument. Specifically,
\begin{align*}
	J(t+1) F_n(u)
	={}& - e^{ix^2/2(t+1)} J (t+1) p e^{-ix^2/2(t+1)} F_n(u) \\
	={}& - e^{ix^2/2(t+1)} J(t+1) p  F_n(e^{-ix^2/2n(t+1)} u) \\
	={}& \frac1n \left(\frac{\partial F_{n}}{\partial z}(u) \left(J\big(n(t+1)\big)u\right) - \frac{\partial F_{n}}{\partial \overline{z}}(u) \left(\overline{J\big(n(t+1)\big)u}\right)\right)
\end{align*}
and, plugging this identity with the Duhamel formula, we have
\begin{align*}
	J(t+1)u(t) ={}& U(t)J(1) u_0 \\
	&{}-\frac{i}n \int_0^t U(t-s) \left( \frac{\partial F_{n}}{\partial z}(u) \left(J\big(n (\cdot+1)\big)u\right) - \frac{\partial F_{n}}{\partial \overline{z}}(u) \left(\overline{J\big(n(\cdot+1)\big)u}\right) \right)(s) ds
\end{align*}
Upon closer inspection, it may seem that obtaining a solution such that $J(t+1)u(t) \in L^2$ by a fixed point argument requires $J(n(t+1))u(t) \in L^2$. However, imposing this requirement leads us to $u_0=u(0) \in X^1(1) \cap X^1(n) = H^{1,1}$.
Hence, when taking data $u_0$ from $X^1(1) \setminus H^{1,1}$, expecting $J(n(t+1))u(t) \in L^2$ becomes impractical.
So, it may seem to be difficult to close estimates with the standard Duhamel formula.

The local well-posedness in $H^{1,1}$ follows from this approach. Indeed, for instance, the local well-posedness in $H^1$ is established in the standard way (see \cite{Ca}).
By substituting
the identity
\[
	J(n(s+1)) u = J(s+1) u - (n-1)(s+1)p u
\]
into the formulation above,
we obtain a closed estimate for $J(t+1)u$.
By a persistence-type argument, we obtain a solution in the $H^{1,1}$-framework.

However, it is important to note that this argument does not guarantee the global-in-time boundedness of $pu$ and $J(t+1)u$. This limitation arises because the identity involves a time-growth factor, $(n-1)(s+1)$. Consequently, we lack a long-time control of solutions even for small data. Thus, achieving global existence and scattering for small data appears challenging with this formulation.

\subsection{A new formulation}
Here, we formulate the equation in a different way than the usual Duhamel formula and obtain a solution such that $J(t+1)u(t) \in L^2$
 as in the gauge-invariant case.

To make the point of the new formulation more comprehensible, we introduce a pseudo-conformal transform to rewrite the problem.
Throughout this subsection, we let $t_0=1$, i.e., $u_0\in X^1(1)$.
We define a variable $v$ as follows: For $t>0$, let
\begin{equation}\label{E:vdef}
	{v}(t,x)  :=
t^{-\frac{d}2} e^{i \frac{|x|^2}{2t}} \overline{ u\left(\frac{1-t}{t}, \frac{x}{t}\right)}.
\end{equation}
Then, as is well-known, if $u$ is a solution to \eqref{1} on an interval $(T_1,T_2)$, where $-1<T_1<T_{2} \leq \infty$, then $v$ is a solution to
\begin{equation}\label{E:v1}
\left\{
\begin{aligned}
&  i \partial _t v =
\frac{p^2}{2}  v + t^{\alpha_0-2} \tilde{F}(t,v), \\
& v(1,x) = e^{i \frac{|x|^2}2} \overline{u_0(x)}
\end{aligned}
\right.
\end{equation}
on $(\frac{1}{T_2+1},\frac{1}{T_1+1}) \subset (0,\infty)$, where
\[
	\tilde{F}(t,v) := \sum_{n \in \mathbb{Z}} \overline{ \lambda _ n} e^{-i \frac{(n-1) |x|^2}{2t} }  F_n  (v).
\]
The standard Duhamel formula of \eqref{E:v1} is
\begin{align}\label{E:Iv1}
	v(t) = U(t-1) v(1) -i \sum_{n =1}^\infty  \overline {\lambda _ n}I_n (t),
\end{align}
where
\begin{equation}\label{E:Indef}
	I_n(t)=I_n(t,v):=  \int_1^t s^{\alpha_0-2} U(t-s) e^{-i \frac{(n-1) |x|^2}{2s} } F_n( v(s) )
 ds .
\end{equation}
Here, the equalities
\[
	\| u(t) \|_{L^2} = \|  v(\tfrac1{1+t})\|_{L^2}
,\quad
	\| J(t+1)u(t) \|_{L^2} = \| \nabla v(\tfrac1{1+t})\|_{L^2}
\]
imply that the solution in the weighted space of \eqref{1} corresponds to an $H^1$-solution of \eqref{E:v1}. Note that the crucial point is that the function space to which the solution belongs no longer changes with time.
Furthermore, the scattering of $u(t)$ as $t\to \infty$ corresponds to the existence of a limit $v(t)$ as $t\to 0+$.

In this formulation, the aforementioned difficulty in the analysis in the weighted space in our setting turns into the following: When the nonlinear term is a gauge-invariant one, i.e., when $\tilde{F}(t,v)= \overline{\lambda _1}F_1(v)$, the analysis in the usual $H^1$ space works well. However, when the nonlinear term includes non-gauge-invariant part $F_n$ ($n\neq1$), it might seem that the presence of quadratic oscillatory terms in the nonlinear term necessitates boundedness of $v$ in a suitable weighted space to evaluate $\nabla \tilde{F}$.
As a result, if we formulate \eqref{E:v1}
in the Duhamel formula,
the Duhamel term does not necessarily make sense as a Bochner integral in $H^1$.

Here we introduce a different formulation  for \eqref{E:Iv1}.
We formally rewrite $I_n$ for $n\ge2$ as follows:
We introduce an operator
\begin{equation}\label{E:Vndef}
	V_n(t) := U(-t)  e^{-i \frac{(n-1) |x|^2}{2t} }.
\end{equation}
It is easy to observe that
\begin{align} \label{3}
i \frac{d}{dt} V_n(t) = -H_n(t) V_n(t)
=- U(-t) \left(\frac{p^2}2 + \frac{n-1}{2t^2}x^2 \right) e^{-i \frac{(n-1) |x|^2}{2t} }
= -V_n(t) \tilde{H}_n(t),
\end{align}
where
\[
	H_n (t):=   \frac{1}{2} p^2 + \frac{(n-1)}{2t^2} \left( x+tp \right)^2 , \quad \tilde{H}_n (t):= \frac{1}2\left(p-\frac{n-1}tx\right)^2 + \frac{(n-1)}{2t^2} x^2 .
\]
The condition $n\ge 2$ implies that the operators $H_n(t)$ and $\tilde{H}_n(t)$ are similar to a harmonic oscillator with a time-dependent coefficients.
Let us further introduce a resolvent operator
\begin{equation}\label{E:Rdef}
 R_n( t) := ( \alpha_0-1+i tH_n(t) ) ^{-1}.
\end{equation}
With the operator, one has
\[
	t^{\alpha_0-2}V_n(t) =  R_n(t) \frac{d}{dt} (t^{\alpha_0-1} V_n(t)).
\]
The identity enables us to do a formal integration by parts:
\begin{align*}
	I_n(t)
 &= U(t) \int_1^t R_n(s) \frac{d}{ds} (s^{\alpha_0-1} V_n(s)) F_n(v(s)) ds\\
&=- U(t) \int_{1}^{t}  \left( \frac{d}{ds}  R_n(s)\right) s^{\alpha_0-1} V_n(s) F_n(v(s)) ds\\
&\quad - U(t) \int_{1}^{t} R_n(s)s^{\alpha_0-1} V_n(s) \left( \frac{d}{ds} F_n(v(s)) \right) ds \\
&\quad + t^{\alpha_0-1} U(t) R_n(t) V_n(t) F_n(v(t)) - U(t) R_n(1)V_n(1)F_n(v(1)).
\end{align*}
Furthermore,
by \eqref{E:v1}, one has
\begin{align*}
	i\frac{d}{ds} F_n(v)
	={}&\frac{\partial F_n}{\partial z} (v)
	\left(\frac{p^2}2 v + s^{\alpha_0-2} \tilde{F}(s,v)\right)
	- \frac{\partial F_n}{\partial \overline{z}} (v) \overline{\left(\frac{p^2}2 v + s^{\alpha_0-2} \tilde{F}(s,v) \right)}.
\end{align*}
According to these identities,
we define
\begin{equation}\label{E:A1def}
	A_{1,n}(t):=- U(t) R_n(1)V_n(1)F_n(v(1)),
\end{equation}
\begin{equation}\label{E:A2def}
	A_{2,n}(t)=A_{2,n}(t,v):=t^{\alpha_0-1} U(t) R_n(t) V_n(t) F_n(v(t)),
\end{equation}
\begin{equation}\label{E:A3def}
	A_{3,n}(t)=A_{3,n}(t,v):=- U(t) \int_{1}^{t} s^{\alpha_0-1} \left( \frac{d}{ds}  R_n(s)\right)  V_n(s) F_n(v(s)) ds,
\end{equation}
\begin{multline}\label{E:A4def}
	A_{4,n}(t)=A_{4,n}(t,v)\\
	:=\frac{i}2 U(t) \int_{1}^{t}  s^{\alpha_0-1}  R_n(s) V_n(s) \left( \frac{\partial F_n}{\partial z} (v(s))
	p^2v(s)- \frac{\partial F_n}{\partial \overline{z}} (v(s))
	\overline{p^2v(s)}\right) ds
\end{multline}
and
\begin{multline}\label{E:A5def}
	A_{5,n}(t)=A_{5,n}(t,v)\\
	:=i U(t) \int_{1}^{t}  s^{2\alpha_0-3}  R_n(s) V_n(s) \left( \frac{\partial F_n}{\partial z} (v(s))
	\tilde{F}(s,v(s))- \frac{\partial F_n}{\partial \overline{z}} (v(s))
	\overline{\tilde{F}(s,v(s))}\right) ds.
\end{multline}
With these functions, we reach to the integral equation
\begin{equation}\label{E:IE}
	v(t) = U(t-1) v(1) - i\lambda_1 I_1(t) -i \sum_{n =2}^\infty \overline{ \lambda _ n} (A_{1,n}(t) + A_{2,n}(t) + A_{3,n}(t)  + A_{4,n}(t) + A_{5,n}(t)).
\end{equation}
We look for a $H^1$-solution to this integral equation (see Proposition \ref{P:key}).

The crucial step of our argument is to derive a closed estimate in the $H^1$-framework by exploiting the properties of the resolvent operator. In terms of the present notation, the aforementioned difficulty can be attributed to the fact that $V_n(s)$ ($n\neq1$) is not a bounded operator on $H^1$. Nevertheless, by using the fact that it is combined with the resolvent operator or the derivative of it, we can recover the boundedness in Sobolev spaces without any weights.
Indeed, we have, for instance,
\begin{align} \label{K2/18-3}
\begin{aligned}
&\left\|  |\sqrt{t}p| U(t) R_n(t) V_n(t) \left| \frac{x}{\sqrt{t}} \right| \right\|_{\CAL{B}(L^r)} \lesssim_r  n^{- \frac14 }, \\
&\left\|  |\sqrt{t}p| U(t) R_n(t) V_n(t) \left\langle \sqrt{t}p \right\rangle \right\|_{\CAL{B}(L^r)} \lesssim_r  n^{\frac12 }
\end{aligned}
\end{align}
for any $1<r<\infty$, $t>0$, and $n\ge2$.
These and similar estimates will be established in Section \ref{ss:bounds}.
These estimates enable us to
close the estimates in the $H^1$-framework (see Lemma \ref{L:Inest}).
Once we obtain a solution to \eqref{E:IE}, we see that it is also a solution to \eqref{E:Iv1}.
The uniqueness needs somewhat different type of argument than usual (see Section \ref{3.3.1}).

In the new formulation,
the problem that the Duhamel term defined by the nonlinear term $t^{\alpha_0-2}\tilde{F}$ does not make sense as a Bochner integral in $H^1$ is resolved as follows: In the new formulation, if $v$ is a continuous function in the $H^1$ sense and belongs to an appropriate auxiliary space, then for all $n \geq 2$, the integrals in $A_{3,n}(t)$, $A_{4,n}(t)$, and $A_{5,n}(t)$  are defined as Bochner integrals in $H^1$.
This means that for $j=3,4,5$, $U(-t)A_{j,n}(t)$ is absolutely continuous as an $H^1$-valued function. Furthermore, $U(-t)A_{1,n}(t)$ is constant with respect to time and is therefore absolutely continuous.
Hence, $U(-t)I_n(t)$ is absolutely continuous in $H^1$ if and only if so is $U(-t)A_{2,n}(t)$.
In summary, the problem
is rephrased as the failure of the absolute continuity of $U(-t)A_{2,n}(t)$.
We note that the failure of the absolute continuity makes the application of Strichartz estimates difficult.
In our proof, we  estimate $A_{2,n}$ without Strichartz estimate.
This is the most crucial point of solving the problem, allowing us to construct an $H^1$ solution that is not absolutely continuous.

On the other hand, a drawback of our new formulation is the appearance of
the term $(\partial_z F_n)(v(s)) p^2 v(s) $ and a similar term, in $A_{4,n}$. The term $(\partial_z F_n)(v(s)) p^2 v$ is one of the most restrictive terms in the evaluation, as it involves the second derivative of $u$. To handle this term, we focus on the resolvent operator, which introduces two antiderivatives, and find that
$R_n(s) V_n(s)(\partial_z F_n)(v(s)) p^2 v$ can be defined in the form sense. Thanks to this fact, we
 handle 
 $(\partial_z F_n)(v(s)) p^2 v$ in the $H^1$-framework
by utilizing the duality argument and the integration by parts.
However, the justification of this procedure produces the limitation $\alpha \geq 2$
because we need to estimate the commutator between $p$ and $(\partial F)(v)$.
We remark that this restriction also produces a restriction on the dimensions.
In fact, our method is applicable to the higher dimension case $d \geq 4$. However, because the restriction $\alpha \geq 2$ is still required, we can handle only mass-critical or -supercritical equations if $d\ge4$, in which case the usual Sobolev space is more suitable to work with than the weighted space in view of scaling.
Moreover, it seems challenging to handle this term in the $H^s$-framework for $s<1$.

Furthermore, since $A_{2,n}$ is a nonlinear term which does not involve any time integration, obtaining an effective small factor from the shortness of the considered time interval is difficult. Therefore, even for \emph{local} well-posedness, small initial data is required. This constitutes the second drawback of our approach.

\bigskip

The rest of the paper is organized as follows. In Section \ref{S:pre}, we introduce estimates associated with the resolvent operator. Additionally, we present basic commutator calculations, decomposition methods using unitary operators, interpolation theorems, and other techniques. In Section \ref{S:proof}, we prove our main theorem.
We first demonstrate that
a solution to the integral equation \eqref{E:IE} is that to \eqref{1}
and then find a solution to \eqref{E:IE} by the fixed point argument
in an appropriate metric space.

\section{Preliminaries}\label{S:pre}

\subsection{Notations}
Let us first collect the notations used throughout this paper.
We denote by  $C$ various positive constants. It is independent of any parameters under consideration.
$A \lesssim B$ implies that $A \le CB$ with some constant $C>0$.
Also $A \lesssim_{r} B$ denotes $A \le C_{r} B$ for a constant $C_{r}$ depending on $r$.
Let $\J{ \cdot} = (1 + | \cdot |^2 )^{1/2}$.
We write $p = -i \nabla _x$. With the notation, one has $- \Delta = |p|^2 = p^2 $. 
We abbreviate
 $L^{r} (\mathbb{R} ^d)$ as $L^r$, for short.
Similarly, $\| \cdot  \|_{L^r}  = \| \cdot  \|_{L^r(\mathbb{R}^d)} $.
We denote $u(t, \cdot)$ as $u(t)$ for short.

The operator norm onto $L^r $, $r \in [1, \infty]$ is denoted as
 $\| \cdot \|_{\CAL{B}(L^r)} $, and inner product on $L^2(\mathbb{R}^d)$ is denoted $(\cdot, \cdot )$.
Let $[ L_1 , L_2 ] =: L_1L_2 -L_2L_1 $ be a commutator between two linear operators $L_1$ and $L_2$.

\subsection{Basic relations on operators}
Our approach employs operator calculation and a useful interpolation theorem (see, e.g., \cite{Ka}). We herein introduce some important relations and lemmas. Introduce a multiplication operator $\CAL{M}(t)$ by
\[
	\CAL{M}(t) u = e^{i \frac{x^2}{2t}} u
\]
for $t\neq0$.
Let us introduce
\[
	A := \frac12 \left( x \cdot p + p \cdot x \right).
\]
Let $\CAL{D}(t)$ be a dilation group such that
\[
( \CAL{D}(t) f )(x) = (it)^{- \frac{d}{2}} f(x/t) .
\]
One has $\CAL{D}(t)^{-1} = i^{d}\CAL{D}(t^{-1})$.
Note that the following representation holds for $t\neq0$:
\begin{align*}
 \CAL{D}(t) f  =  e^{-\frac{t d\pi}{4|t|}i } \left( e^{-i (\log |t|) A } f \right) 
\end{align*}
(see, e.g.,  Isozaki \cite{Iso}).
The following relations are standard:
\begin{align*}
i[ p_j , x_k ] =  \delta _{jk}, \quad i[ A, x_j]= x_j, \quad i[A, p_j ] = - p_j.
\end{align*}
The following identities are frequently used:
\begin{equation} \label{7}
\begin{aligned}
 e^{-ia  x^2} \MAT{x \\ p } e^{ia x^2} &= \MAT{ x \\ p + 2a x },&
e^{-ia  p^2} \MAT{x \\ p } e^{ia p^2} &= \MAT{ x-2ap \\ p }, \\
 e^{i a  A} \MAT{x \\ p } e^{-ia  A} &= \MAT{ e^{a} x \\ e^{-a} p },&
\mathcal{D}(a)^{-1} \MAT{x \\ p } \mathcal{D}(a) &= \MAT{ a x \\ a^{-1} p }.
\end{aligned}
\end{equation}
The proof of \eqref{7} is simple. For instance, it is done by looking at the derivative in $a$ for each terms. Indeed, the third identities is proved as follows:
For $\phi \in \SCR{S}({\bf R}^d)$,
\begin{align*}
\frac{d}{da}e^{i a  A} \MAT{x_j \\ p_j } e^{-ia  A} \phi = e^{i a  A} \MAT{ i [ A, x_j] \\ i[A, p_j] } e^{-ia  A} \phi = e^{i a  A} \MAT{ x_j \\ -p_j } e^{-ia  A} \phi
\end{align*}
holds.
This equation is an ordinary differential equation for the vector quantity ${\displaystyle e^{i a  A} \MAT{x_j \\ p_j } e^{-ia  A}  \phi }$. Solving them, we obtain the desired identity from $e^{i 0  A}  x_j e ^{-i 0  A} \phi = x_j \phi $ and $e^{i 0  A}  p_j e ^{-i 0  A} \phi = p_j \phi $.
The first two identities are established in a similar way.
The last one is an immediate consequence of the third.
\begin{Rem}
One quick application of \eqref{7} is the calculation of the
 commutation between $e^{-ia p^2} $ and $  x \cdot p $:
\begin{align*}
e^{-ia p^2}  x \cdot p = (e^{-ia p^2} x e^{ia p^2} )\cdot (e^{-ia p^2} p e^{ia p^2})   e^{-ia p^2} = (x-2ap) \cdot p e^{-ia p^2} .
\end{align*}
\end{Rem}

\subsection{$L^r$ bound for resolvent operators}
In our argument, the resolvent operator
\begin{align*}
R_{n, \mathrm{os}} := \left( \alpha_0 -1 +i \sqrt{n-1} H_{\mathrm{os}} /2 \right) ^{-1}, \quad H_{\mathrm{os}} = p^2 + x^2
\end{align*}
plays very important role, where $n\ge2$.
We herein introduce the $L^r$-bounds of several operators associated with $R_{n, \mathrm{os}}$.
\begin{Lem} \label{L2/16}
For all $r \in (1, \infty)$, we have the bounds
\begin{align} \label{K2/15-2}
\begin{aligned}
	& \left\| p^2 R_{n, \mathrm{os}}  \right\|_{\mathcal{B}(L^{r})} \lesssim n^{- \frac12}, \quad \left\| x^2 R_{n, \mathrm{os}}  \right\|_{\mathcal{B}(L^{r})} \lesssim n^{- \frac12},  \quad \left\| (x \cdot p + p \cdot x) R_{n, \mathrm{os}}  \right\|_{\mathcal{B}(L^{r})} \lesssim n^{- \frac12}, \\
	& \left\|  R_{n, \mathrm{os}}  p^2 \right\|_{\mathcal{B}(L^{r})} \lesssim n^{- \frac12}, \quad \left\|  R_{n, \mathrm{os}} x^2 \right\|_{\mathcal{B}(L^{r})} \lesssim n^{- \frac12},  \quad \left\| R_{n, \mathrm{os}} (x \cdot p + p \cdot x)  \right\|_{\mathcal{B}(L^{r})} \lesssim n^{- \frac12}
\end{aligned}
\end{align}
and
\begin{align} \label{K2/16-1}
 \left\| R_{n, \mathrm{os}}^{} \right\|_{\CAL{B}(L^r)} \lesssim n^{- \frac12}
\end{align}
for $n\ge2$.
\end{Lem}
\begin{proof}[Proof of Lemma \ref{L2/16}]
Noting that all the following symbols
\begin{align*}
|x|^2 (\xi ^2 + x^2 +i)^{-1}, \quad    | \xi |^2 (\xi ^2 + x^2 +i)^{-1}, \quad (x \cdot \xi + \xi \cdot x) (\xi ^2 + x^2 +i)^{-1}
\end{align*}
belong to the H\"{o}lmander class  $S_{1,0}^0 $, we see from Taylor \cite{Tay} (see the comment below Proposition 5.5, of \cite{Tay} and see, also Cardona-Delgado-Ruzhanski \cite{CDR}) that
for all $1< r<\infty$,
\begin{align*}
p^2 (p^2+x^2 +i)^{-1}, \ x^2 (p^2+x^2 +i)^{-1} , \  (x \cdot p + p \cdot x)(p^2+x^2 +i)^{-1} \in \CAL{B}(L^r).
\end{align*}
The latter three terms in \eqref{K2/15-2} can be handled as well. We show \eqref{K2/16-1}. We first let $2 \leq n \leq N_0$ for some $N_0 \gg 1$ to be chosen later. Then by the same argument in the proof of \eqref{K2/15-2}, we find \eqref{K2/16-1} for $n \leq N_0$ as $\| R_{n, \mathrm{os}} \| _{\CAL{B}(L^r)} \leq C_{N_0} n^{- \frac12}$. Next we consider the case where $n \geq N_0 \gg 1$.  Note that $H_{\mathrm{os}}^{-1} $ is bounded operator on $L^r({\mathbb{R} ^d})$ for any $r \in [1, \infty]$ (see, \cite[Theorem 1]{BT}) and whose bound is independent of $n$. Hence one finds by the second resolvent identity that
\begin{align*}
R_{n , \mathrm{os}} -(i \sqrt{n-1} H_{\mathrm{os}}/2)^{-1} &=
 R_{n, \mathrm{os}} \left( i\sqrt{n-1} H_{\mathrm{os}}/2 - R_{n, \mathrm{os}}^{-1}  \right)  (i \sqrt{n-1} H_{\mathrm{os}}/2)^{-1}
\\ & = (1- \alpha_0) R_{n, \mathrm{os}} ( \sqrt{n-1} H_{\mathrm{os}}/2)^{-1}.
\end{align*}
By taking $N_0$ large enough so that $2(\alpha_0 -1) \| H_{\mathrm{os}}^{-1}  \|_{\CAL{B}(L^r)} (N_0 -1)^{-\frac12} < 1/2  $, we have \eqref{K2/16-1}.
\end{proof}
For $a,b \in \mathbb{R}$ and $0 \leq \theta \leq 2$, define an operator
\begin{align*}
| ax + b p |^{\theta} := |a|^{\theta} e^{ibp^2/(2a)} |x|^{\theta} e^{-ib p^2/(2a)} = |b|^{\theta} e^{-iax^2/(2b)} |p|^{\theta} e^{ia x^2/(2b)} .
\end{align*}
Moreover, on $ C_0^{\infty} (\mathbb{R}^d) \subset L^2(\mathbb{R}^d)$, we define the $R_{n, \mathrm{os}}^{\theta}$ as the pseudo-differential operator $R_{n, \mathrm{os}}^{\theta} (x,p)$ with symbol
\begin{align*}
 r_{n , \mathrm{os}} (x, \xi) =  \left( \alpha_0 -1 -i \sqrt{n-1}  \left( \xi ^2 + x^2 \right)  /2 \right)^{-\theta} .
\end{align*}

Then, we have the following bounds:
\begin{Lem} \label{L4}
For all $r \in (1, \infty)$ and $0 \leq \theta \leq 1$, we have the bounds
\begin{align} \label{10/30-1}
\begin{aligned}
 \left\| |p|^{2\theta}  R_{n, \mathrm{os}}^{\theta}  \right\|_{\mathcal{B}(L^{r})} \lesssim{}& n^{- \frac{\theta}2}, \\
 \left\|  R_{n, \mathrm{os}}^{\theta}  |p+ \sqrt{n-1} x|^{2\theta} \right\|_{\mathcal{B}(L^{r})} \lesssim{}& n^{ \frac{\theta}2}, \\
 \left\| |x|^{2\theta} R_{n, \mathrm{os}}^{\theta}  \right\|_{\mathcal{B}(L^{r})} + \left\|  R_{n, \mathrm{os}}^{\theta} |x|^{2\theta} \right\|_{\mathcal{B}(L^{r})} \lesssim{}& n^{- \frac{\theta}2}
\end{aligned}
\end{align}
for $n\ge2$.
\end{Lem}

This lemma immediately follows from the following interpolation theorem:
\begin{Lem}[Interpolation theorem] \label{L3}
Let $r \in (1,\infty)$. Let $n\ge2$ and let
$P=|ax +bp|^2$ and $Q=R_{n, \mathrm{os}}$ or $P= R_{n, \mathrm{os}} $ and $  Q=|ax +bp|^2$. Then, for all $\theta \in [0,1]$,
\begin{align} \label{K2/15-1}
\left\| P^{\theta}  Q^{\theta}   \right\|_{\CAL{B}(L^r)}  \lesssim  \left( \frac{a^2 + b^2}{\sqrt{n}} \right) ^{\theta}
\end{align}
holds. Moreover,
\begin{align} \label{K2/17-1}
\| R_{n, \mathrm{os}}^{\theta} \|_{\CAL{B}(L^r)} \lesssim n^{- \frac{\theta}2}
\end{align}
for $n\ge2$.
\end{Lem}
\Proof{ We only show the case where
$P=|ax +bp|^2$ and $Q=R_{n, \mathrm{os}} $ with $b \neq 0$. Other cases can be shown similarly. We know $Q^{\theta} \in \CAL{B}(L^r)$ since $(\alpha_0 -1 + \sqrt{n-1}(|x|^2 + |\xi|^2  )/2 )^{- \theta} \in S_{1,0}^{-2 \theta} $, and $P^{\theta}= |b|^{2\theta} e^{-i \frac{a |x|^2 }{2b}}  |p|^{2\theta} e^{i \frac{a |x|^2 }{2b}} $.
Here we define $$P_{\rho}^{\theta} =  |b|^{2\theta} e^{-iax^2/(2b)} (\chi (|p| / \rho) |p|)^{2\theta} e^{ia x^2/(2b)}$$
for $\rho >0$, where $\chi \in C_0^{\infty} ({\mathbb R})$. Then we find that $P_{\rho}^{\theta} \in \CAL{B}(L^r)$ by Young's convolution inequality, and $P_{\rho} \to P$, $\rho \to \infty$ on the suitable function space. Hence, to show \eqref{K2/15-1}, it suffices to show the uniform bound of
\begin{align*}
  \left\| (\chi (|p| / \rho ) |bp|)^{2\theta} \cdot e^{ia x^2/(2b)} \cdot Q^{\theta}  \right\| _{\CAL{B}(L^r)}. 
\end{align*}
One obtains the bound
by interpolation as in \cite{Ka} and \cite{Iso} (see also \cite{DKP}).
If we obtain the following two bounds
\begin{align*}
  \left\| (\chi (|p| / \rho ) |bp|)^{0} \cdot e^{ia x^2/(2b)} \cdot Q^{0}  \right\| _{\CAL{B}(L^r)} &\leq C_0, \\   \left\| ( \chi (|p| / \rho ) |bp|)^{2} \cdot e^{ia x^2/(2b)} \cdot Q^{1}  \right\| _{\CAL{B}(L^r)} &\leq C_1,
\end{align*}
then $C_0^{1- \theta} C_1^{\theta} $ is the desired uniform bound.
Clearly the first bound holds with the choice $C_0 =1$. Further, by using $b^2 p^2 e^{ia x^2/(2b)}=e^{ia x^2/(2b)}(a^2 x^2 + b^2 p^2 + ab (x \cdot p + p \cdot x )) $ and  \eqref{K2/15-2}, we have the second estimate with $C_1 = C(a^2+b^2) n^{- \frac12}$.
The estimate \eqref{K2/17-1} can be shown similarly, thanks to \eqref{K2/16-1}.
}

\subsection{Estimates on resolvent operators}\label{ss:bounds}
Let us collect estimates on the resolvent operator.
We introduce one more operator
\begin{equation}\label{E:tRdef}
	\tilde{R}(t):= V_n(t)^{-1} R_n(t) V_n(t),
\end{equation}
where $R_n(t) = (\alpha_0 -1 + i t H_n(t))^{-1} $, $H_n(t) = \frac{p^2}2 + \frac{(n-1)}{2t^2} \left( x+ tp \right)^2  $, $\alpha _0 = \frac{d}2(\alpha -1) >1$ and $V_n(t) = U(-t) e^{-i \frac{(n-1) |x| ^2}{2t}}$. Due to the relation
$\tilde{H}_n(t)=V_n(t)^{-1}H_n(t) V_n(t)
$, which follows from \eqref{3}, one sees that
$\tilde{R}(t)$ is the resolvent operator of $\tilde{H}_n(t)$, i.e.,
\begin{equation}\label{E:tRdef2}
 \tilde{R}_n( t) = ( \alpha_0-1+i t\tilde{H}_n(t) ) ^{-1}.
\end{equation}

Let us collect estimates on the resolvent operators
 $R_n(t)$ and $\tilde{R}_n(t)$, here we recall that
 \begin{align*}
 R_n(t) = (\alpha_0 -1 + it H_n(t) )^{-1}
 &= \left(  \alpha_0 -1 +  \frac{it}2 \left( p^2 + (n-1)\frac{(x+ tp)^2 }{t^2} \right)  \right) ^{-1}
 \\ &= U(-t) \left(  \alpha_0 -1 +  \frac{it}2 \left( p^2 + (n-1)\frac{x^2 }{t^2} \right)  \right) ^{-1} U(t).
 \end{align*}
To obtain bounds on $R_n(t)$ and $\tilde{R}_n(t)$, we initially introduce the following factorization of the resolvent operator:

\begin{Lem} \label{L:Rest}
Suppose that $n \geq 2$. For any fixed $t >0$, $R_n(t)$ and $\tilde{R}_n(t)$  given in \eqref{E:Rdef} and \eqref{E:tRdef}, respectively, have the following representation:
\begin{align}\label{E:Rnform}
R_n(t)  &= U(-t)\mathcal{D}((\tfrac{t^2}{n-1})^{\frac14})
R_{n, \mathrm{os}}
\mathcal{D}((\tfrac{t^2}{n-1})^{-\frac14}) U(t)
\end{align}
and
\begin{align}\label{E:tRnform}
\tilde{R}_n(t)  &= e^{ i\frac{(n-1)|x|^2}{2t} }\mathcal{D}((\tfrac{t^2}{n-1})^{\frac14})R_{n, \mathrm{os}}
\mathcal{D}((\tfrac{t^2}{n-1})^{-\frac14})
e^{- i\frac{(n-1)|x|^2}{2t} }.
\end{align}
\end{Lem}
\begin{proof}
By \eqref{7} (or \eqref{3}), we find that
\begin{align*}
H_n(t) &= \frac12 U(-t) \left(p^2 + \frac{n-1}{t^2}x^2 \right) U(t)
\\ &= \frac{\sqrt{n-1}}{2t} U(-t)
\mathcal{D}((\tfrac{t^2}{n-1})^{\frac14})
H_{\mathrm{os}}
\mathcal{D}((\tfrac{t^2}{n-1})^{-\frac14})
U(t).
\end{align*}
The formula \eqref{E:Rnform} easily follows.
\eqref{E:tRnform} is similarly valid and thus the proof is completed.
\end{proof}

Based on the formula \eqref{E:Rnform} and \eqref{E:tRnform}, we introduce the powers of $R_n(t)$ and $\tilde{R}_n(t)$, respectively.
For $\gamma \in \mathbb{R}$, we let
\[
	R_n(t)^\gamma := U(-t)\mathcal{D}((\tfrac{t^2}{n-1})^{\frac14})
R_{n, \mathrm{os}}^{\gamma}
\mathcal{D}((\tfrac{t^2}{n-1})^{-\frac14}) U(t)
\]
and
\[
	\tilde{R}_n(t)^\gamma := e^{ i\frac{(n-1)|x|^2}{2t} }\mathcal{D}((\tfrac{t^2}{n-1})^{\frac14})R_{n, \mathrm{os}}^{\gamma}
\mathcal{D}((\tfrac{t^2}{n-1})^{-\frac14})
e^{- i\frac{(n-1)|x|^2}{2t} }.
\]

The crucial step of the proof is to use the resolvent to replace derivative by weight and to gain the derivative at the sacrifice of time decay factor.
The following lemma, which is a consequence of Lemma \ref{L:Rest}, is one of the main tool in our argument.
\begin{Lem} \label{L:Rest2}
Fix $t>0$. The following estimates hold for $\gamma\in [0,2]$, $r\in(1,\infty)$, and $n\ge2$:
\begin{align} \label{E:Rest3}
	\left\| | \sqrt{t} p |^{\gamma} U(t) R_n(t)^{\frac\gamma2} U(-t)
\right\|_{\mathcal{B}(L^r)} \lesssim{}& 1, \\
	\left\| \left|\frac{ x}{\sqrt{t}} \right|^{\gamma} U(t)R_n(t)^{\frac\gamma2} V_n(t) \right\|_{\CAL{B}(L^r)} \lesssim{}& n^{ -\frac{\gamma}2}. \label{E:Rest4}
\end{align}
Further, following estimates are valid:
\begin{align}
 \label{E:Rest0} \left\| U(t)R_n(t)^{\frac\gamma2} V_n(t)    \right\|_{\CAL{B}(L^r)} = \left\|    \tilde{R}_n(t)  ^{\frac\gamma2}    \right\| _{\CAL{B}(L^r)}
	&\lesssim n^{ -\frac{\gamma}4},  \\
 \label{8} \left\| U(t)R_n(t)^{\frac\gamma2} V_n(t) \left|\frac{ x}{\sqrt{t}} \right|^{\gamma}     \right\|_{\CAL{B}(L^r)} = \left\|    \tilde{R}_n(t)  ^{\frac\gamma2}  \left|\frac{ x}{\sqrt{t}} \right|^{\gamma}   \right\| _{\CAL{B}(L^r)}
	&\lesssim n^{ -\frac{\gamma}2},  \\
 \label{11} \left\| U(t) R_n(t)^{\frac\gamma2} V_n(t) |\sqrt{t} p|^{\gamma }    \right\| _{\CAL{B}(L^r)}  = \left\|  \tilde{R}_n(s)  ^{ \frac\gamma2} | \sqrt{t}p|^{\gamma }  \right\| _{\CAL{B}(L^r)}
	&\lesssim n^{ \frac{\gamma}2}.
\end{align}
\end{Lem}
\begin{proof}
Let us begin with the proof of \eqref{E:Rest3}. By
\begin{align*}
\left\| | \sqrt{t} p |^{\gamma} U(t) R_n(t)^{\frac\gamma2} U(-t)
\right\|_{\mathcal{B}(L^r)} = \left\| \left|  {p} { (n-1)^{\frac14} } \right|^{\gamma} R_{n, \mathrm{os}}^{\frac\gamma2}
\right\|_{\mathcal{B}(L^r)} \lesssim 1,
\end{align*}
where we have used \eqref{10/30-1}.
Similarly, \eqref{E:Rest4} and \eqref{E:Rest0} follows from \eqref{10/30-1} and \eqref{K2/17-1}, respectively.
Let us next obtain the bound in \eqref{8}.
By \eqref{E:tRdef},
\eqref{E:tRnform}, and \eqref{7}, one deduces that
\begin{align*}
	U(t) R_n(t)^{\frac{\gamma}2} V_n(t)  \left|\frac{x}{\sqrt{t}}\right|^\gamma
	={}& e^{ -i\frac{(n-1)|x|^2}{2t} } \tilde{R}_n(t) ^{\frac{\gamma}2} \left|\frac{x}{\sqrt{t}}\right|^\gamma \\
	={}& \mathcal{D}((\tfrac{t^2}{n-1})^{\frac14}) R_{n, \mathrm{os}}^{\frac{\gamma}{2}} \left|\frac{x}{{(n-1)}^{1/4}}\right|^\gamma
	\mathcal{D}((\tfrac{t^2}{n-1})^{-\frac14})e^{- i\frac{(n-1)|x|^2}{2t} }.
\end{align*}
The first identity shows
the first identity of \eqref{8}, and the second inequality of \eqref{8} follows from $\eqref{10/30-1}$.
The bound in \eqref{11} is obtained in a similar way.
We only note that
\begin{align*}
	U(t) R_n(t)^{\frac{\gamma}{2}} V_n(t) |\sqrt{t}p|^{\gamma}
	={}&  e^{ -i\frac{(n-1)|x|^2}{2t} }\tilde{R}_n(t)^{\frac{\gamma}{2}} |\sqrt{t}p|^{\gamma} \\
	={}& \mathcal{D}((\tfrac{t^2}{n-1})^{\frac14})R_{n, \mathrm{os}}^{\frac{\gamma}{2}}
{(n-1)^{\frac{\gamma}4}}\left| p+\sqrt{n-1} x \right|^{\gamma}\mathcal{D}((\tfrac{t^2}{n-1})^{-\frac14})e^{- i\frac{(n-1)|x|^2}{2t} }
\end{align*}
holds.
\end{proof}

\begin{Rem}
The first estimate of \eqref{K2/18-3} follows from \eqref{E:Rest3}, \eqref{E:Rest0}, and \eqref{8} with the choice $\gamma=1$.
Similarly, \eqref{E:Rest3}, \eqref{E:Rest0}, and \eqref{11}
give us the second estimate.
\end{Rem}
The next assertion comes from the duality of Lemma \ref{L:Rest2}.
\begin{Lem} \label{L:Rest2d}
Fix $t >0$ and $n \ge 2$.
Let $r\in(1,\infty)$ and $ \gamma \in [0,2]$. Then the following estimates hold:
\begin{align}
	\norm{U(t) \overline{R_n(t)^{\frac\gamma2}} U(-t) |\sqrt{t} p |^{\gamma}}_{\mathcal{B}(L^r)} \lesssim{}& 1, \label{E:Rest3d} \\
	\norm{V_n(t)^{-1}\overline{R_n(t)^{\frac{\gamma}{2}}} U(-t)\abs{\frac{ x}{\sqrt{t}}}^{\gamma}}_{\mathcal{B}(L^{r})} \lesssim{}& n^{-\frac{\gamma}{2}}, \label{E:Rest4d}
\end{align}
Moreover these estimates
\begin{align}
	\norm{V_n(t)^{-1} \overline{R_n(t)^{\frac{\gamma}{2}}} U(-t)}_{\mathcal{B}(L^{r})} \lesssim n^{-\frac{\gamma}{4}}, \label{E:Rest0d} \\
	\norm{\abs{\sqrt{t}p}^{\gamma} V_n(t)^{-1}\overline{R_n(t)^{\frac{\gamma}{2}}} U(-t)}_{\mathcal{B}(L^{r})} \lesssim n^{\frac{\gamma}{2}} \label{11d}
\end{align}
are valid.
\end{Lem}

We also need to handle the derivative of $R_n(t)$.
\begin{Lem} \label{L-K5/4-1}
Let $\theta_1, \theta_2 \geq 0$ satisfy $\theta_1 + \theta_2  \leq 2$. Let $r\in (1,\infty)$.
It holds that
\begin{align} \label{E:Rest2}
	\left\| U(t) R_n(t)  ^{-\frac{\theta_1}2} \left( \frac{d}{dt}  R_n(t)  \right) R_n(t) ^{-\frac{{\theta_2}}2}  U(-t)\right\|_{\mathcal{B}(L^r)}
	\lesssim n^{\frac{{\theta_1}+{\theta_2} +2}4} t^{-1}
\end{align}
for all $t\in (0,1]$ and $n\ge2$.
\end{Lem}
\begin{proof}
We denote
$
	P_n(t) = U(-t)\mathcal{D}((\tfrac{t^2}{n-1})^{\frac14}).
$
Using the identities
\begin{align*}
i \frac{d}{dt} P_n(t)^{-1} &=
	-\frac{1}{2t}\left( A  - \sqrt{ n-1} p^2  \right) P_n(t)^{-1}
\end{align*}
and
\begin{align*}
	i \frac{d}{dt} P_n(t) &=
	\frac{1}{2t} P_n(t) \left( A  - \sqrt{ n-1}p^2  \right) ,
\end{align*}
one obtains
\begin{align*}
	i\frac{d}{ds} R_n(t)
	={}&  \frac{1}{2s} P_n(t) \left[ A - \sqrt{ n-1}p^2 , (\alpha_0-1+i \tfrac{\sqrt{n-1}}{2}H_{\mathrm{os}} ) ^{-1} \right]
P_n(t)^{-1}\\
 ={}&  -\frac{1}{2t}  P_n(t) R_{n, \mathrm{os}}\left[ A - \sqrt{ n-1}p^2 , \alpha_0-1+i\tfrac{\sqrt{n-1}}{2} H_{\mathrm{os}} \right]
R_{n, \mathrm{os}}P_n(t)^{-1}.
\end{align*}
Due to the identities
\begin{align*}
i[A, p^2]= -2 p^2, \quad i[A, x^2] = 2 x^2, \quad i[p^2, p^2] =0, \quad i[p^2 , x^2] = 4 A,
\end{align*}
one obtains
\[
	\left[  A  + \sqrt{ n-1}p^2 , \alpha_0-1+i \tfrac{\sqrt{n-1}}{2}H_{\mathrm{os}} \right] = -\sqrt{n-1} (p^2-x^2) + 2(n-1)A.
\]
Thus, if $\theta_1,\theta_2\ge0$ satisfy $\theta_1+\theta_2\le2$ then
we have
\begin{align*}
& \left\| U(t) R_n(t)  ^{-\frac{\theta_1}2} \left( \frac{d}{dt}  R_n(t)  \right) R_n(t) ^{-\frac{{\theta_2}}2}  U(-t)\right\|_{\mathcal{B}(L^r)}
\\ & \leq  C\J{n}^{\frac12}t^{-1}   \left\| R_{n, \mathrm{os}}^{1 - {\theta_1}/2}  \left( -p^2 + x^2 - 2\sqrt{n-1} A \right) R_{n, \mathrm{os}}^{ 1 - {\theta_2} /2}   \right\|_{\mathcal{B}(L^r)}.
\end{align*}
Hence, the desired bound follows if we have the bound
\[
 \left\| R_{n, \mathrm{os}}^{1 - {\theta_1}/2}  \left( -p^2 + x^2 - 2\sqrt{n-1} A \right) R_{n, \mathrm{os}}^{ 1 - {\theta_2} /2}   \right\|_{\mathcal{B}(L^r)} \lesssim n^{\frac{{\theta_1} + {\theta_2}}4}.
\]
Let us show the bound.
We consider the case $1- {\theta_1} /2 < 1/2$ and $1- {\theta_2} /2 > 1/2$, for instance.
The other case is handled in a similar way.
Using the identity
\[
 - p^2 +x^2 -2 \sqrt{n-1}A = - (p+ \sqrt{n-1} x)^2 + n x^2,
\]
we have
\begin{align*}
& \left\| R_{n, \mathrm{os}}^{1 - {\theta_1}/2}  \left( -p^2 + x^2 - 2\sqrt{n-1} A \right) R_{n, \mathrm{os}}^{ 1 - {\theta_2} /2}   \right\|_{\mathcal{B}(L^r)}
\\ & \lesssim \sup_{f\in L^r, g \in L^{r'}} \Big\{ \left\| \left| p + \sqrt{n-1}x \right|  ^{{\theta_1}} R_{n, \mathrm{os}}^{1- \frac{{\theta_2}}2} f \right\|_{L^r} \left\| \left| p + \sqrt{n-1}x\right| ^{2-{\theta_1}} R_{n, \mathrm{os}}^{1- \frac{{\theta_1}}2} g  \right\|_{L^{r'}}
 \\
& \qquad + n \left\| |x|^{{\theta_1}}  R_{n, \mathrm{os}}^{1- \frac{{\theta_2}}2} f \right\|_{L^r}    \left\|  |x | ^{2-{\theta_1}} R_{n, \mathrm{os}}^{1- \frac{{\theta_1}}2} g  \right\|_{L^{r'}} \Big\}
\\ & \lesssim n^{\frac{{\theta_1}}4 - \frac12 + \frac{{\theta_2}}4 + \frac{{\theta_1}}4 + \frac12 - \frac{{\theta_1}}4} + n^{1- \frac{{\theta_1}}4 - \frac12 + \frac{{\theta_2}}4 + \frac{{\theta_1}}4 - \frac12 + \frac{{\theta_1}}4} =2 n^{\frac{{\theta_1}+ {\theta_2}}4},
\end{align*}
where we have used \eqref{10/30-1} and \eqref{K2/17-1}.
This completes the proof of \eqref{E:Rest2}.
\end{proof}

\subsection{The pseudo-conformal transform}

Let us give a notion of solution to \eqref{E:v1}.
\begin{Def}\label{D:v1sol}
Let $I \subset (0,\infty)$ be a interval containing 1.
We say a function $v(t,x)$ is a solution to \eqref{E:v1} on $I$ if $v \in C(I;H^{1}) \cap L^{q_0}_{\mathrm{loc}}(I;L^{r_0})$ and the identity 
\begin{align}\label{E:Iv2}
	v(t) = U(t-1) v(1) -i \int_1^t s^{\alpha_0-2} U(t-s) \tilde{F}(s,v(s))ds
\end{align}
holds for $t\in I$.
\end{Def}
We remark that if $v \in C(I;H^{1})$
then, for all $n$, $I_n$ defined in \eqref{E:Indef} makes sense as a Bochner integral in $L^2$. (One sees that $I_1$ makes sense also an integral in $H^1$.) It will turn out in our proof that $I_n$ belongs
to $C(I;H^1)$ for any $n\ge2$.
Further, together with the weighted summability assumption, we see that
$\sum_{n=2}^\infty \overline{\lambda_n} I_n$ belongs to $C(I;H^1)$.
As a result, we will see that if $v$ is a solution on $I$ in the above sense then the identity
 \eqref{E:Iv1} also holds on $I$.

\begin{Lem}\label{L7}
Fix $T>0$.
$u(t,x)$ is a solution to \eqref{1} on an interval $[0,T)$ in the sense of Definition \ref{D:sol} if and only if
 $v(t,x)$ given in \eqref{E:vdef}
is a solution to \eqref{E:v1}
on $(\frac1{T+1},1]$
in the sense of Definition \ref{D:v1sol}.
Further, the following identities hold for $t\in [0,T)$:
\begin{align*}
\| u(t) \|_{L^{2}} &{}= \| v(\tfrac1{t+1}) \|_{L^{2}}, &
\| J(t+1) u(t)\|_{L^2}
&{}=
 \| p v(\tfrac1{t+1}) \|_{L^2}  .
\end{align*}
Moreover, $u(t)$ scatters to $U(-t-1)u_+$
in $L^2$ as $t\to\infty$ if and only if
$\lim_{t\to0+} v(t)$ exists and equals to  $i^{\frac{d}2} \mathcal{F}^{-1} \overline{u_+}$
in $L^2$.
Furthermore,
$u(t)$ scatters
as $t\to \infty$ in the sense of Theorem \ref{T1}
 if and only if
$\lim_{t\to0+}v(t)$ exists in the $H^{1}$-topology.
\end{Lem}
\Proof{
It is well-known (see, e.g., \cite{Ca}).
Let us give a brief proof, for self-containedness.
One sees from \eqref{E:vdef} that
the following identity holds:
\[
	U(-t-1) u(t) = i^{\frac{d}2} e^{-i \frac{|x|^2}{2(t+1)}} \mathcal{F}^{-1} \overline{v(\tfrac1{t+1})}.
\]
Hence, one deduces that $\|u(t)\|_{L^2} = \|v(\tfrac1{t+1})\|_{L^2}$ and that
$U(-t-1) u(t)\to u_+$ in $L^2$ as $t\to \infty$ if and only if $v(t) \to i^{\frac{d}2} \mathcal{F}^{-1} \overline{u_+}$ in $L^2$.
Similarly,
\[
	U(-t-1)J(t
+1) u(t)
	=xU(-t-1) u(t)
	= i^{\frac{d}2} e^{-i \frac{|x|^2}{2(t+1)}} \mathcal{F}^{-1} \overline{(p v)(\tfrac1{t+1})}
\]
follows. Hence, one has
$\| J(t+1) u(t)\|_{L^2} =
\| p v(\tfrac1{t+1}) \|_{L^2} $
and the desired statement on the equivalence between the scattering of $u(t)$ as $t\to\infty$ and the existence of a limit of $v(t)$ as $t\to0+$.
}

\subsection{Fractional Hardy inequality}

The following estimate is well-known (\cite{He77}, see also \cite{RuSu}).
\begin{Lem}
Let $d\ge1$ and
$r \in (1,\infty)$. If $0 \le s < \frac{d}r$ then
\[
	\| |x|^{-s} f\|_{L^r(\mathbb{R}^d)} \lesssim
	\| |p|^s f \|_{L^r(\mathbb{R}^d)}
\]
for any $f \in H^s_r(\mathbb{R}^d)$.
\end{Lem}

\section{Proof of the main theorem}\label{S:proof}

In this section, we prove our main result.
In Sections \ref{ss:formulation} to \ref{ss:keyprop}, we prove the GWP part of Theorem \ref{T1}.
Section \ref{ss:scattering} is devoted to the proof of the scattering.
We prove the scattering-part of Theorem \ref{T1} and the small data scattering in $H^{1,1}$ for both time directions (Corollary \ref{C:both}).

\subsection{Formulation}\label{ss:formulation}

Thanks to the argument in Remark \ref{R:red}, we restrict ourselves to the case $t_0=1$ in the proof of Theorem \ref{T1}.
We use the standard fixed point argument to find a solution to \eqref{E:IE}.
We regard the right hand side of \eqref{E:IE} as a functional $\Phi(v)$, i.e.,
\begin{equation}\label{E:Phidef}
\begin{aligned}
	\Phi(v)&:= U(t-1) v(1) - i {\overline{\lambda_1}} I_1(t,v) \\&\quad -i \sum_{n =2}^\infty  {\overline{\lambda_n}} (A_{1,n}(t) + A_{2,n}(t,v) + A_{3,n}(t,v)  + A_{4,n}(t,v) + A_{5,n}(t,v)),
\end{aligned}
\end{equation}
where we regard $v(1)=e^{i \frac{|x|^2}2} \overline{u_0} \in H^1$, which appears in the first term and $A_{1,n}$, as a fixed function.
For a function $v$, we denote
\begin{equation}\label{E:error}
\CAL{E}(v)=\CAL{E}(v)(t,x):= i\partial_t v - \frac{p^2}2 v - t^{\alpha_0-2}\tilde{F}(t,v).
\end{equation}
By definition of $\Phi$, one has
\begin{equation}\label{E:Phierror}
\begin{aligned}
	&{}i \partial_t \Phi(v) - \frac{p^2}2 \Phi(v) \\
	= {}& t^{\alpha_0-2} \tilde{F}(t,v)
	+ \sum_{n\ge2} {\overline{\lambda_n}}  t^{\alpha_0-1} U(t) R_n(t) V_n(t) \left( \frac{\partial F_n}{\partial z}(v) \CAL{E}(v) - \frac{\partial F_n}{\partial \overline{z}}(v) \overline{\CAL{E}(v)} \right).
\end{aligned}
\end{equation}

We define a complete metric space as follows:
For $M>0$, we define
\[
	Z_M:= \{ v \in C((0,1],H^1) ; \|v\|_X + \|v\|_Y \le M \},
\]
where
\[
	\| v \|_{X} :=  \| v \|_{L^\infty ((0,1];H^{1})} + \|v\|_{L^{q_0} ((0,1];H^{1}_{r_0})} , \quad
	\| v \|_{Y} := \|\partial_t v\|_{L^{q_1} ((0,1];H^{-1}_{r_1})}
\]
with
\[
	(q_0,r_0) = \begin{cases}
	(12,3) & d=1,\\
	(\frac{2}{\alpha_0-1},\frac2{2-\alpha_0}) & d=2, \\
	(4, 3) & d=3.
	\end{cases}
\]
and
\[
	( q_1 , r_1  ) :=
\begin{cases}
(1, 2) & d=1, \\
(1, r_0) & d=2, 3.
\end{cases}
\]
The metric is defined as
\[
	d_Z(v_1,v_2) := \|v_1-v_2\|_{X}.
\]
Obviously, $(Z_M,d_Z)$ is a complete metric space for all $M\ge0$.
By virtue of Lemma \ref{L7},
the global-existence part of Theorem \ref{T1} follows from the following proposition:
\begin{Prop}\label{P:key}
Under the assumption of Theorem \ref{T1}, there exists $\ep_0$ and $M>0$
such that the following assertions hold:
\begin{enumerate}
\item
If $u_0 \in L^2$ satisfies
$u_0 \in J (1) L^2$ and
$\|u_0\|_{L^2} + \| J(1) u_0 \|_{ L^2 } \le \ep_0 $ then
$\Phi$ is a contraction map from $(Z_M,d_Z)$ to itself.
\item If $v \in Z_M$ is a fixed point of
$\Phi$ then $\CAL{E}(v)=0$, i.e., $v$ is a solution to \eqref{E:v1} on $(0,1]$.
\end{enumerate}
\end{Prop}

\subsection{Contraction property}
In this subsection, we prove the first assertion of Proposition \ref{P:key}. Namely,
we shall show that, under a suitable choice of $\ep_0$ and $M$, the map $\Phi$ given in
\eqref{E:Phidef} is a contraction map from $(Z_M,d_Z)$ to itself.
To this end, we show that
\[
	\|\Phi(v)\|_{X}  + \|\Phi(v)\|_{Y} \le M
\]
for all $v\in Z_M$ and
\[
	d_Z(\Phi(v_1),\Phi(v_2)) \le \frac12 d_Z(v_1,v_2)
\]
for all $v_1,v_2\in Z_M$.

\subsubsection{The relation between $X$-norm and $Y$-norm}
We first claim that the estimates on the $Y$-norm essentially follow from those on the $X$-norm.
To this end, we obtain an estimate on the error term.
\begin{Lem} \label{lem:32}
Fix $\ep_0>0$ and $M>0$. Then,
\[
	\|\CAL{E}(v)\|_{L^{q_1} ((0,1];H^{-1}_{r_1})} \lesssim M(1+M^{\alpha-1})
\]
for $v\in Z_M$.
\end{Lem}
\begin{proof}
Pick $v\in Z_M$.
$i \partial_t v \in L^{q_1} ((0,1]; H^{-1}_{r_1})$ follows from $\|v\|_Y\le M$. Further,
\[
	\| \tfrac{p^2}2 v \|_{L^{q_1} ((0,1]; H^{-1}_{r_1})}
	\lesssim \| v \|_{L^{q_1} ((0,1];H^{1}_{r_1})}  \le \| v \|_X \le M
\]
holds from the embedding $L^{q} ((0,1]) \hookrightarrow L^{1} ((0,1])$ ($\forall q \in [1,\infty]$).
Finally, the embedding $ L^\infty ((0,1],H^1) \hookrightarrow L^\infty ((0,1],L^{\alpha r_1})$ yields
\begin{align*}
	\| t^{\alpha_0-2} \tilde{F} (t,v) \|_{L^1((0,1];H^{-1}_{r_1})}
	&\lesssim \| t^{\alpha_0-2} \|_{L^1((0,1])}
	\| |v|^\alpha \|_{L^\infty((0,1];L^{r_1})} \\
	&\lesssim 	\| v \|_{L^\infty((0,1];H^1)}^\alpha \\ & \lesssim M^\alpha .
\end{align*}
Thus, we obtain the desired estimate.
\end{proof}

\begin{Lem}\label{L:Phierror}
Fix $\ep_0>0$ and $M>0$. Then,
\[
	\left\|i \partial_t \Phi(v) - \frac{p^2}2 \Phi(v)\right\|_{L^{q_1} ((0,1];H^{-1}_{r_1})} \lesssim M^{\alpha}+ M^{\alpha-1} \|\CAL{E}(v)\|_{L^{q_1} ((0,1];H^{-1}_{r_1})}
\]
for $v\in Z_M$.
\end{Lem}
\begin{proof}
Recall the identity \eqref{E:Phierror}. As in the proof of Lemma \ref{lem:32}, one has
\[
	\| t^{\alpha_0-2} \tilde{F} (t,v) \|_{L^1((0,1];H^{-1}_{r_1})}
	\lesssim M^\alpha.
\]
Let us estimate the second term of the right hand side of \eqref{E:Phierror}.
By duality,
\[
	\left\| U(t) R_n(t)V_n(t)
	\frac{\partial F_n}{\partial z} (v)  \CAL{E} \right\|_{H^{-1}_{r_1}}
	\lesssim  \sup_{\| f \|_{H^1_{r_1'}} =1 }   \left\| \overline{\frac{\partial F_n}{\partial z} (v)} V_n(t)^{-1} \overline{R_n(t)} U(-t)   f \right\|_{H^1_{r_{1}'}}  \|\CAL{E}\|_{H^{-1}_{r_1}}.
\]
It follows that
\begin{align}\label{K2/18-2}
\begin{aligned}
	\left\| \frac{\partial F_n}{\partial z} (v) V_n(t)^{-1} \overline{R_n(t)} U(-t)   f \right\|_{H^1_{r_1'}}
	\lesssim{}& \left\| \left( \nabla \frac{\partial F_n}{\partial z} (v) \right)  V_n(t)^{-1} \overline{R_n(t)} U(-t)   f \right\|_{L^{r_1'}} \\
&{}+ \left\| \frac{\partial F_n}{\partial z} (v)  \J{p}V_n(t)^{-1} \overline{R_n(t)} U(-t)   f \right\|_{L^{r_1'}} .
\end{aligned}
\end{align}

Let us consider the case $d=1$.
Recall that $r_1=2$.
Thanks to \eqref{11d},
the first term of the right hand side of \eqref{K2/18-2} is bounded from above by a constant multiple of
\begin{align*}
	\| v \|_{L^{\infty}}^{\alpha -2}  \| v \|_{H^1} \left\| |p| V_n(t)^{-1} \overline{R_n(t)} U(-t)   f  \right\|_{L^{2}}^{\frac12}  \left\|  V_n(t)^{-1} \overline{R_n(t)} U(-t)   f \right\|_{L^2}^{\frac12} \lesssim t^{ - \frac{1}4} M^{\alpha -1}.
\end{align*}
Similarly, combining \eqref{11d} with \eqref{E:Rest4d} and \eqref{E:Rest0d}, the second term is bounded by
\begin{align*}
 t^{ -\frac12} \left\| v\right\|_{L^{\infty}}^{\alpha -1}  \left\| \J{\frac{x}{\sqrt{t} }}^{-1} f \right\|_{L^2} \lesssim  t^{-\frac14} \left\| v\right\|_{L^{\infty}}^{\alpha -1}  \left\| f \right\|_{L^{\infty}} \lesssim
t^{-\frac{1}4} M^{\alpha -1}
\end{align*}
with a constant.
Thus, it follows that
\[
	\left\|  t^{\alpha_0-1} U(t) R_n(t) V_n(t)  \frac{\partial F_n}{\partial z} \CAL{E}(v)   \right\|_{L^1((0,1];H^{-1})} \le CM^{\alpha-1} \| \CAL{E} \|_{L^1((0,1];H^{-1})}.
\]
The term with $\frac{\partial F_n}{\partial \overline{z}} \overline{\CAL{E}(v)}$ is handled similarly.
Thus, summing up in $n$, we obtain the result, where we have used $\sum n^2 |\lambda_n| < \infty$.

Let us consider the case $d=2$. By means of the H\"{o}lder inequality, 
the first term of the right hand side of \eqref{K2/18-2} is bounded by
\[
\left\|  v\right\|_{H^1} \left\| v \right\|_{L^4}^{\alpha -2} \left\|  V_n(t)^{-1} \overline{R_n(t)} U(-t) f \right\|_{L^{\frac{4}{\alpha_0 -1}}}
\lesssim M^{\alpha-1}\left\| U(t) \overline{R_n(t)} U(-t) f \right\|_{L^{\frac{4}{\alpha_0 -1}}}.
\]
By \eqref{E:Rest3d}, \eqref{E:Rest0d} and Sobolev embedding, one has
\[
	\left\|  U(t) \overline{R_n(t)} U(-t) f \right\|_{L^{\frac{4}{\alpha_0 -1}}}
	\lesssim \| |\sqrt{t}p|^{-\frac{\alpha_{0}-1}{2}} f \|_{L^{\frac{4}{\alpha_0 -1}}} \lesssim t^{-\frac{\alpha_0-1}4}
	\| f\|_{H^{1}_{r_0'}}
	\le t^{-\frac{\alpha_0-1}4}.
\]
By using \eqref{11d} with \eqref{E:Rest4d} and \eqref{E:Rest0d}, one sees that the second term is bounded by
\begin{align*}
t^{ - \frac12 +  \frac{3-2 \alpha_0}{2}}  \left\| v \right\|_{L^{\frac{\alpha_0}{\alpha_0 -1}}}^{\alpha-1} \left\| |x|^{  -(3-2 \alpha_0)  }  f \right\|_{L^{\frac{2}{2-\alpha_0}}}  \lesssim t^{1-\alpha_0} M^{\alpha -1} \| f \|_{H^1_{r_0'}}
\le t^{1-\alpha_0} M^{\alpha -1}.
\end{align*}
Thus,
\[
	\left\|  t^{\alpha_0-1} U(t) R_n(t) V_n(t)  \frac{\partial F_n}{\partial z} \CAL{E}(v) \right\|_{L^{1}((0,1);H^{-1}_{r_0})} \lesssim M^{\alpha-1} \| \CAL{E} \|_{L^{1}((0,1);H^{-1}_{r_0})}.
\]
The term with $ \frac{\partial F_n}{\partial \overline{z}} \overline{\CAL{E}(v)}$ is handled similarly.

Let us consider the case $d=3$. As the same calculations in the case of $d=2$, the first term of the right hand side of \eqref{K2/18-2} is bounded by
\[
\left\|  v\right\|_{H^1}\left\|  V_n(t)^{-1} \overline{R_n(t)} U(-t) f \right\|_{L^{6}}
\lesssim M^{}\left\| U(t) \overline{R_n(t)} U(-t) f \right\|_{L^{6}} \lesssim Mt^{- \frac14} \left\| |p|^{- \frac12} f \right\|_{L^{6}}  \lesssim t^{- \frac14}M .
\]
The second term is bounded by
\begin{align*}
	t^{- \frac12} \left\| v\right\|_{L^6} \left\|f \right\|_{L^2}
	\lesssim  t^{- \frac12} M.
\end{align*}
Thus,
\[
	\left\|  t^{\alpha_0-1} U(t) R_n(t) V_n(t)  \frac{\partial F_n}{\partial z} \CAL{E}(v) \right\|_{L^{1}((0,1);H^{-1}_{r_0})} \lesssim M^{\alpha-1} \| \CAL{E} \|_{L^{1}((0,1);H^{-1}_{r_0})},
\]
where we remark that $\alpha _0 -1 = \frac12$ when $d=3$.

\end{proof}

Utilizing these two lemmas,
we see that
\begin{equation}\label{E:XtoY}
\begin{aligned}
	\| \Phi(v)
	\|_{Y} &\le \frac12
	\left\| p^2 \Phi(v) \right\|_{L^{q_1}((0,1];H^{-1}_{r_1})} + \left\| i\partial_t\Phi(v) -\frac{p^2}2 \Phi(v) \right\|_{L^{q_1}((0,1];H^{-1}_{r_1})}
	\\
	&\lesssim  \left\|  \Phi(v)\right\|_{X} + M^\alpha(1+M^{\alpha-1}).
\end{aligned}
\end{equation}
Note that we have used the embedding
$L^q_t ((0,1]) \hookrightarrow L^1_t ((0,1])$ for all $q\in[1,\infty]$ to estimate the first term.
This inequality means that the estimate on $Y$-norm of $\Phi(v)$ follows from that on $X$-norm.

\subsubsection{Key estimates}
We turn to the estimate of $\| \Phi(v) \|_X$ and $d_Z(\Phi(v_1),\Phi(v_2))$.
Let us begin the proof with the estimate of the linear part and $I_1$ part of $\Phi$.
They are handled by a standard argument.
By Strichartz and
Lemma \ref{L7}, one has
\begin{equation}\label{E:Pkeypf1}
	\|U(t-1)v(1)\|_X \lesssim
	\|v(1)\|_{H^1} \lesssim \|u_0\|_{L^2} + \| J(1)u_0\|_{L^2 } \lesssim \ep_0.
\end{equation}
This term does not appear in the estimate of $\Phi(v_1)-\Phi(v_2)$ since this part is independent of the choice of $v_j \in Z_M$.

The $X$-norm of $I_1$ is estimated as follows:
If $d=1$ then, by using Strichartz and the fact that $H^1$ is an algebra, one has
\begin{equation}\label{E:Pkeypf2}
\begin{aligned}
	\| I_1(\cdot,v) \|_{X}
	\lesssim \|s^{\alpha_0-2} F_1(v)\|_{L^{1}_s((0,1]; H^{1})}
\lesssim \| s^{\alpha_0-2} \|_{L^{1}_s((0,1])}
	\|v\|_{L^\infty((0,1];H^{1})}^{\alpha}\lesssim \|v\|_{X}^\alpha \le M^\alpha
\end{aligned}
\end{equation}
for $v\in Z_M$.
Further,
\begin{equation}\label{E:Pkeypf2.5}
\begin{aligned}
	&\|I_1(\cdot,v_1)-I_{1}(\cdot,v_2)\|_{X} \\
&\lesssim
	\|s^{\alpha_0-2} (F_1(v_1)-F_1(v_2))\|_{L^{1}_s((0,1]; H^1)} \\
& \lesssim  (\|v_1\|_{L^\infty((0,1]\times \mathbb{R})}+\|v_2\|_{L^\infty((0,1]\times \mathbb{R})})^{\alpha-1} \|v_1-v_2\|_X\\
& \quad + (\|v_1\|_{L^\infty((0,1]\times \mathbb{R})}+\|v_2\|_{L^\infty((0,1]\times \mathbb{R})})^{\alpha-2} (\|v_1\|_X+\|v_2\|_X) \|v_1-v_2\|_{L^\infty((0,1]\times \mathbb{R})}\\
&\lesssim M^{\alpha-1} \|v_1-v_2\|_X
\end{aligned}
\end{equation}
for $v_1,v_2 \in Z_M$.

If $d=2$, $3$ then we see from Strichartz and the embedding $H^1 \hookrightarrow L^p$ for $p\in [2,\infty)$ when $d=2$ and  $p\in [2,6]$ when $d=3$ that
\begin{equation}\label{E:Pkeypf3}
\begin{aligned}
	\| I_1(\cdot,v) \|_{X}
	\lesssim{}& \|s^{\alpha_0-2} F_1(v)\|_{L^{q_0'}_s((0,1]; H^{1}_{r_0'})} \\
\lesssim{}& \| s^{\alpha_0-2} \|_{L^{\frac{2}{3-\alpha_0}}_s((0,1])}
	\|v\|_{L^\infty((0,1];L^{\frac{2\alpha_0}{\alpha_0-1}})}^{\alpha-1} \|v\|_{L^\infty ((0,1];H^1)} \\
\lesssim{}& \|v\|_{X}^\alpha \le M^\alpha
\end{aligned}
\end{equation}
for $v\in Z_M$ and
\begin{equation}\label{E:Pkeypf4}
\begin{aligned}
	&\|I_1(\cdot,v_1)-I_{1}(\cdot,v_2)\|_{X}\\
	&\lesssim \|s^{\alpha_0-2} (F_1(v_1)-F_1(v_2))\|_{L^{q_0'}_s((0,1]; H^1_{r_0'})} \\
&\lesssim (\|v_1\|_{L^\infty((0,1];L^{\frac{2\alpha_0}{\alpha_0-1}})}+\|v_2\|_{L^\infty((0,1];L^{\frac{2\alpha_0}{\alpha_0-1}})})^{\alpha-1} \|v_1-v_2\|_X \\
&\quad+ (\|v_1\|_{L^\infty((0,1];L^{\frac{2\alpha_0}{\alpha_0-1}})}+\|v_2\|_{L^\infty((0,1];L^{\frac{2\alpha_0}{\alpha_0-1}})})^{\alpha-2}(\|v_1\|_X+\|v_2\|) \|v_1-v_2\|_{L^\infty((0,1];L^{\frac{2\alpha_0}{\alpha_0-1}})} \\
&\lesssim M^{\alpha-1} \|v_1-v_2\|_X
\end{aligned}
\end{equation}
for $v_1,v_2 \in Z_M$.

It will turn out that we have a similar estimate for $\left\| I_n \right\|_{X}$ for $n\ge2$.
\begin{Lem}\label{L:Inest}
For $n\ge2$, one has
\begin{equation}\label{E:Inest1}
 \sum_{k=1}^5\| A_{k,n}(\cdot,v) \|_{X}  \lesssim n^\frac52(\ep_0^{\alpha} + M^{\alpha} \J{M}^{\alpha-1})
\end{equation}
for $v\in Z_M$
and
\begin{equation}\label{E:Inest2}
	\sum_{k=1}^5\| A_{k,n}(\cdot,v_1)-A_{k,n}(\cdot,v_2)\|_X  \lesssim n^\frac72M^{\alpha-1}\J{M}^{\alpha-1} \|v_1-v_2\|_X
\end{equation}
for $v_1,v_2 \in Z_M$.
\end{Lem}
We prove this lemma in the rest of this subsection.

\subsubsection{Proof of Lemma \ref{L:Inest} in one dimension}

We shall prove Lemma \ref{L:Inest} in the one dimensional case.

Recall that $(q_0,r_0) = \left(12, 3\right)$.
Fix $n\ge2$ and pick
$v,v_1,v_2\in Z_M$.
Let us first estimate $A_{1,n}(t)$ defined in \eqref{E:A1def}.
Recall that this term is independent of $v\in Z_M$.
By Strichartz estimate,
\[
	\| A_{1,n} \|_{X}
\lesssim \| R_n(1) V_n(1) F_n(v(1))\|_{H^{1}}.
\]
One sees from \eqref{E:Rest3} and \eqref{E:Rest0} that
\[
	\left\| \J{p} R_n(1) V_n(1)\right\|_{L^{2}}
	=\left\| \J{p} U(1)R_n(1) U(-1)\right\|_{L^{2}}
	\lesssim 1.
\]
Hence,
\begin{equation}\label{E:d1A1n1}
	\| R_n(1) V_n(1) F_n(v(1))\|_{H^{1}}
	\lesssim
\|(\CAL{M}(1) u(0))^\alpha \|_{L^2} \lesssim  \|u(0)\|_{L^2}^{\frac{\alpha(\alpha_0+1)}{2\alpha_0+1}}
\|J(1)u(0)\|_{L^2}^{\frac{\alpha\alpha_0}{2\alpha_0+1}} \lesssim \ep_0^\alpha.
\end{equation}

Let us consider $A_{2,n}(t,v)$ 
given in \eqref{E:A2def}.
Note that
\[
\|\J{p} \J{t^{1/2}p}^{-1}\|_{\mathcal{B}(L^{r})} \lesssim t^{-\frac12}
\]
and
\[
	\left\| \J{\sqrt{t} p} U(t) R_n(t)^{\frac{1}2}U(-t) \cdot U(t) R_n(t)^{\frac{1}2} V_n(t) \J{\frac{ x}{\sqrt{t}}} \right\|_{\mathcal{B}(L^{r})} \lesssim 1
\]
for $r=2,3$.
One has
\begin{align*}
	\| A_{2,n}(t) \|_{H^1_r}
	\lesssim t^{\alpha_0-\frac32}\left\| \J{\frac{ x}{\sqrt{t}}}^{-1}  \right\|_{L^r}
\| F_n(v(t))\|_{L^\infty}
	= C t^{\alpha_0-\frac32 + \frac1{2r}}
\| v(t)\|_{L^\infty}^\alpha
\end{align*}
for $r=2,3$.
Thus, thanks to the embedding $H^1\hookrightarrow L^\infty$,
we obtain
\begin{equation}\label{E:d1A2n1}
	\| A_{2,n}(\cdot,v) \|_X \lesssim \|v\|_{L^\infty((0,1]\times \mathbb{R})}^\alpha \le M^\alpha
\end{equation}
for $\alpha_0>5/4$.
Similarly,  one has
\begin{equation}\label{E:d1A2n2}
\begin{aligned}
	\|A_{2,n}(\cdot,v_1)-A_{2,n}(\cdot,v_2)\|_X &\lesssim \| F_n(v_1(t))- F_n(v_2(t))\|_{L^\infty((0,1]\times \mathbb{R})}\\
&\lesssim n M^{\alpha-1}\|v_1-v_2\|_X.
\end{aligned}
\end{equation}

Let us turn to the estimate of $A_{3,n}$ given in \eqref{E:A3def}.
By Strichartz,
\[
	\|A_{3,n}(\cdot,v)\|_X \lesssim
 \left\| U(s)
	\left( \frac{d}{ds}  R_n(s)\right) s^{\alpha_0-1} V_n(s) F_n(v(s))
\right\|_{L^{1}((0,1];H^{1})}.
\]
By means of Lemma \ref{L:Rest2} and Lemma \ref{L-K5/4-1},
we have the following boundedness of the operators:
For $s\in (0,1]$,
\[
\|\J{p} \J{s^{1/2}p}^{-1}\|_{\mathcal{B}(L^{2})} \lesssim s^{-\frac12},
\]
\[
	\|\J{\sqrt{s}p} U(s)R_n(s)^{\frac12} U(-s)\|_{\mathcal{B}(L^2)} \lesssim 1,
\]
\[
	\left\| U(s) R_n(s)^{-\frac12}
	\left( \frac{d}{ds}  R_n(s)\right) R_n(s)^{-\frac12}U(-s)\right\|_{\mathcal{B}(L^{2})} \lesssim n^{} s^{-1}
\]
and
\[
	\left\|U(s) R_n(s)^{\frac12} V_n(s)\J{\frac{x}{\sqrt{s}}}
\right\|_{\mathcal{B}(L^{2})} \lesssim 1
.
\]
Therefore, for $s\in(0,1]$,
\begin{align*}
	 \left\| U(s)
	\left( \frac{d}{ds}  R_n(s)\right) s^{\alpha_0-1} V_n(s) F_n(w(s))
\right\|_{H^1}
	&\lesssim n s^{\alpha_0-\frac52} \left\| \J{\frac{x}{\sqrt{s}}}^{-1}\right\|_{L^2} \|   F_n(v(s))  \|_{L^{\infty}}\\
	&=C n s^{\alpha_0- \frac{9}{4}}
	\|v(s)\|^\alpha_{L^\infty} \\
	&\lesssim n s^{\alpha_0- \frac{9}{4}}
	\|v\|_{L^\infty((0,1];H^{1})}^\alpha.
\end{align*}
Since $s^{\alpha_0-\frac{9}{4}} \in L^{1}(0,1)$ for $\alpha_0> \frac{5}{4}$, we have the bound
\begin{equation}\label{E:d1A3n1}
	\|A_{3,n}(\cdot,v)\|_X \lesssim n^{}M^\alpha.
\end{equation}
Further, a similar argument shows
\begin{equation}\label{E:d1A3n2}
\begin{aligned}
	\|A_{3,n}(\cdot,v_1)-A_{3,n}(\cdot,v_2)\|_X &\lesssim n^{}
 \| F_n(v_1(t))- F_n(v_2(t))\|_{L^\infty((0,1]\times \mathbb{R})} \\
&\lesssim n^{2}  M^{\alpha-1} \|v_1-v_2\|_X.
\end{aligned}
\end{equation}

We move on to the estimate of $A_{4,n}(t)$ defined in \eqref{E:A4def}.
By Strichartz,
\begin{align*}
	\|A_{4,n}(\cdot,v)\|_X &\lesssim
 \left\|
	 U(s) R_n(s) s^{\alpha_0-1} V_n(s) \frac{\partial F_n}{\partial z} (v(s)) p^2 v(s)
\right\|_{L^1 ((0,1];H^{1})}\\
&\quad +  \left\|
	 U(s) R_n(s) s^{\alpha_0-1} V_n(s) \frac{\partial F_n}{\partial \overline{z}} (v(s)) \overline{p^2 v(s)}
\right\|_{L^1 ((0,1];H^{1})}.
\end{align*}
We only consider the first term since the second term is handled similarly.
One sees from Lemma \ref{L:Rest2} that
\[
	\left\| \J{\sqrt{s}p} U(s)R_n(s) ^{\frac12} U(-s)\cdot U(s)  R_n(s)^{\frac{1}2} V_n(s) \J{\sqrt{s} p}
 \right\|_{\mathcal{B}(L^{2})} \lesssim n^{\frac12},
\]
which yields
\[
	\left\| U(s)
	 R_n(s) s^{\alpha_0-1} V_n(s) \frac{\partial F_n}{\partial z} (v(s))
	p^2 v(s)
\right\|_{H^1}
\lesssim n^{\frac12} s^{\alpha_0-2} 	\left\|   \J{p}^{-1} \frac{\partial F_n}{\partial z} (v(s)) p^2 v(s) \right\| _{L^2}.
\]
We claim that
\begin{equation}\label{E:d1A4nclaim}
	\left\|   \J{p}^{-1} \frac{\partial F_n}{\partial z} (v(s)) p^2 v(s) \right\| _{L^2}
\lesssim n^2 \|v\|_{H^1}^{\alpha-1}
\| v\|_{H^{1}_{3 }  }.
\end{equation}
Indeed, by duality,
\begin{align*}
		\left\|   \J{p}^{-1} \frac{\partial F_n}{\partial z} (v(s)) p^2 v(s) \right\| _{L^2}
	&=  \sup_{\|\varphi\|_{L^2}=1}
	\left| \left(  p v,  p \left( \overline{\frac{\partial F_n}{\partial z} (v)} \J{p}^{-1}\varphi \right) \right)\right| \\
& \lesssim \|v\|_{H^{1}_3}
\sup_{\|\varphi\|_{L^2}=1}
\left\| p\left( \overline{\frac{\partial F_n}{\partial z} (v)} \J{p}^{-1}\varphi \right) \right\|_{L^{\frac3{2}}}.
\end{align*}
By H\"older and Sobolev,
\begin{align*}
	\sup_{\|\varphi\|_{L^2}=1}
\left\|   \overline{\frac{\partial F_n}{\partial z} (v)} (p\J{p}^{-1}\varphi)  \right\|_{L^{\frac32}} &{}\lesssim n \left\| v \right\|_{L^{12\alpha_0}}^{\alpha-1}  \lesssim n\|v\|_{H^1}^{\alpha-1}.
\end{align*}
By using an estimate for the product of functions in Triebel-Lizorkin space \cite[\S 4.4.4 Theorem 3]{RS}, one obtains
\[
	\sup_{\|\varphi\|_{L^2}=1}\left\|   \overline{\left(p\frac{\partial F_n}{\partial z} (v)\right)} \J{p}^{-1}\varphi  \right\|_{L^{\frac32}}
\lesssim \sup_{\|\varphi\|_{L^2}=1}\left\| p\frac{\partial F_n}{\partial z} (v) \right\|_{L^2}
\|\J{p}^{-1} \varphi\|_{H^1}
\lesssim \left\| p\frac{\partial F_n}{\partial z} (v) \right\|_{L^2}
.
\]
Further, using the embedding $H^1 \hookrightarrow L^\infty$, we obtain
\[
	\left\| p\frac{\partial F_n}{\partial z} (v) \right\|_{L^2}
\lesssim n^2  \|v\|_{H^1}^{\alpha-1}. \]
Thus, we obtain the claim \eqref{E:d1A4nclaim}. Hence,
\begin{equation*}
	\|A_{4,n}(\cdot,v)\|_{X} \lesssim
n^{\frac52} \| s^{\alpha_0-2} \|_{L^{\frac{12}{11}}((0,1])} \|v\|_X^{\alpha-1}.
\end{equation*}
Since $s^{\alpha_0-2} \in L^{\frac{12}{11}}((0,1])$ $\Leftrightarrow$ $\alpha_0>\frac{13}{12}$, one arrives at the estimate
\begin{equation}\label{E:d1A4n1}
		\|A_{4,n}(\cdot,v)\|_{X} \lesssim
n^{\frac52} M^\alpha.
\end{equation}

Similarly, one sees from Strichartz that
\begin{align*}
	&\|A_{4,n}(\cdot,v_1)-A_{4,n}(\cdot,v_2)\|_X \\
&\lesssim n^\frac12
 \left\|
	  s^{\alpha_0-2} \J{p}^{-1} \left( \frac{\partial F_n}{\partial z} (v_1(s)) p^2 v_1(s) - \frac{\partial F_n}{\partial z} (v_2(s)) p^2 v_2(s) \right)
\right\|_{L^1 ((0,1]; L^2)}\\
&\quad + n^\frac12 \left\|
	 s^{\alpha_0-2} \J{p}^{-1} \left(\frac{\partial F_n}{\partial \overline{z}} (v_1(s)) \overline{p^2 v_1(s)}-\frac{\partial F_n}{\partial \overline{z}} (v_2(s)) \overline{p^2 v_2(s)}\right)
\right\|_{L^1 ((0,1]; L^2)}.
\end{align*}
Let us only consider the first term. Mimicking the proof of \eqref{E:d1A4nclaim}, one has
\begin{equation*}
	\left\|   \J{p}^{-1} \frac{\partial F_n}{\partial z} (v_1) p^2 (v_1-v_2) \right\| _{L^2}
\lesssim n^2 \|v_1\|_{H^1}^{\alpha-1}
\| v_1-v_2\|_{H^{1}_{3 }  } \lesssim n^2 M^{\alpha-1} \|v_1-v_2\|_X
\end{equation*}
and
\begin{align*}
	&\left\|   \J{p}^{-1} \left( \frac{\partial F_n}{\partial z} (v_1(s)) -\frac{\partial F_n}{\partial z} (v_2(s)) \right) p^2 v_2(s) \right\| _{L^2}\\
&\lesssim
\| v_2\|_{H^{1}_{3 }  } \left(
\left\|   {\frac{\partial F_n}{\partial z} (v_1)} - {\frac{\partial F_n}{\partial z} (v_2)} \right\|_{L^{6}}
+ \left\| p\left(\frac{\partial F_n}{\partial z} (v_1)-\frac{\partial F_n}{\partial z} (v_2)\right) \right\|_{L^2}
\right) \\
&\lesssim n^3 M^{\alpha-1}\|v_1-v_2\|_X
.
\end{align*}
Hence, one has
\begin{equation}\label{E:d1A4n2}
	\|A_{4,n}(\cdot,v_1)-A_{4,n}(\cdot,v_2)\|_X \lesssim n^{\frac72} M^{\alpha-1} \|v_1-v_2\|_X.
\end{equation}

We finally treat $A_{5,n}(t)$ defined in \eqref{E:A5def}.
By Strichartz,
\begin{align*}
	\|A_{5,n}(\cdot,v)\|_X &\lesssim
 \left\| U(s)
	 R_n(s) s^{2\alpha_0-3} V_n(s) \frac{\partial F_n}{\partial z} (v(s)) \tilde{F}(s,v(s)) \right\|_{L^{1} ((0,1];H^1)}\\
&\quad + \left\| U(s)
	 R_n(s) s^{2\alpha_0-3} V_n(s) \frac{\partial F_n}{\partial \overline{z}} (v(s)) \overline{\tilde{F}(s,v(s))} \right\|_{L^{1} ((0,1];H^1)}.
\end{align*}
We only handle the first term. Estimating as in the proof of \eqref{E:d1A2n1}, one has
\begin{align*}
	\left\| U(s)
	 R_n(s) s^{2\alpha_0-3} V_n(s) \frac{\partial F_n}{\partial z} (v(s)) \tilde{F}(s,v(s)) \right\|_{H^1}
	\lesssim ns^{2\alpha_0-3- \frac12+ \frac{1}4} \|v(s)\|_{L^\infty}^{2\alpha-1}.
\end{align*}
Hence,
\begin{equation}\label{E:d1A5n1}
	\| A_{5,n} \|_X \lesssim n\|v\|_{L^\infty ((0,1]\times\mathbb{R})}^{2\alpha-1} \lesssim nM^{2\alpha-1}
\end{equation}
holds if $\alpha_0>9/8$.
In a similar way, one sees that
\begin{equation}\label{E:d1A5n2}
\begin{aligned}
	\|A_{5,n}(\cdot,v_1)-A_{5,n}(\cdot,v_2)\|_X &\lesssim
 \left\| \frac{\partial F_n}{\partial z} (v_1) \tilde{F}(s,v_1) -\frac{\partial F_n}{\partial z} (v_2) \tilde{F}(s,v_2)\right\|_{L^\infty ((0,1]\times\mathbb{R})}\\
&\quad + \left\| \frac{\partial F_n}{\partial \overline{z}} (v_1) \overline{\tilde{F}(s,v_1)}-\frac{\partial F_n}{\partial \overline{z}} (v_2) \overline{\tilde{F}(s,v_2)} \right\|_{L^\infty ((0,1]\times\mathbb{R})}\\
&\lesssim n^2 M^{2(\alpha-1)} \|v_1-v_2\|_X
\end{aligned}
\end{equation}
if $\alpha_0>9/8$.

The estimate
\eqref{E:Inest1}
follows from
\eqref{E:d1A1n1}, \eqref{E:d1A2n1},
\eqref{E:d1A3n1},
\eqref{E:d1A4n1},
and \eqref{E:d1A5n1}.
Similarly, \eqref{E:d1A2n2},
\eqref{E:d1A3n2},
\eqref{E:d1A4n2},
and \eqref{E:d1A5n2} yield \eqref{E:Inest2}.

\subsubsection{Proof of Lemma \ref{L:Inest} in two dimensions}
We next consider the two dimensional case of the first assertion of Lemma \ref{L:Inest}.
Recall that
\[
	(q_0,r_0) = \left(\frac2{\alpha_0-1}, \frac2{2-\alpha_0}\right).
\]
Fix $n\ge2$ and pick
$v,v_1,v_2\in Z_M$.
Let us first estimate $A_{1,n}(t)$ given in \eqref{E:A1def}.
By Strichartz,
\[
	\| A_{1} \|_{X}
\lesssim \| R_n(1) V_n(1) F_n(v(1))\|_{H^1}.
\]
We see from \eqref{E:Rest3} and \eqref{E:Rest0} that
\[
	\left\| \J{p} R_n(1) V_n(1)\right\|_{L^{2}}
	=\left\| \J{p} U(1)R_n(1) U(-1)\right\|_{L^{2}}
	\lesssim 1,
\]
which yields
\begin{equation}\label{E:d2A1n1}
	\| R_n(1) V_n(1) F_n(v(1))\|_{H^1}
	\lesssim
\| (\CAL{M}(1)u(0))^\alpha \|_{L^2} \lesssim \|u(0)\|_{L^2}^{\frac{\alpha}{\alpha_0+1}}
\|J(1)u(0)\|_{L^2}^{\frac{\alpha \alpha_0}{\alpha_0+1}}
 \le \ep_0^\alpha.
\end{equation}

Let us estimate $A_{2,n}$ given in \eqref{E:A2def}. Since
\[
	\left\| \J{\sqrt{t}p} U(t) R_n(t)^{\frac12}U(-t) \cdot
	U(t) R_n(t)^{\frac{1}{2}} V_n(t) \J{\frac{ x}{\sqrt{t}}} \right\|_{\mathcal{B}(L^{\frac{2}{2-\alpha_0}})} \lesssim 1
\]
follows from \eqref{E:Rest3} and \eqref{8}, one has
\begin{align*}
	\| A_{2,n}(t) \|_{H^1_{\frac{2}{2-\alpha_0}}}
	&{}\lesssim t^{\alpha_0-1} \left\|  \J{p} \J{\sqrt{t} p}^{-1}  \J{\frac{x}{\sqrt{t}}}^{-1} F_n(v(t)) \right\| _{L^{\frac{2}{2-\alpha_0}}} \\
 &{}\lesssim t^{\alpha_0-\frac{3}2+\frac{a}2} \| |x|^{-\frac{a}{\alpha_0+1}} v(t) \|^{\alpha}_{L^{\frac{2 (\alpha_0+1)}{2-\alpha_0}}},
\end{align*}
for $a\in [0,1] $,
where we have used the estimate $\J{
x}^2 \gtrsim |x|^{\theta}$ ($\theta \in [0,2]$).
For $\alpha_0 \in [\frac75,2]$,
we choose $a=0$.
Then, since $t^{\alpha_0-\frac32} \in L^{\frac{2}{\alpha_0-1}}_t((0,1])$ for $\alpha_0>4/3$ and $H^1(\mathbb{R}^2)\hookrightarrow L^p(\mathbb{R}^2)$ for $p\in[2,\infty)$, one obtains
\begin{equation}\label{E:d2A2npf1}
	\| A_{2,n}(\cdot,v) \|_{L^{\frac{2}{\alpha_0-1}}((0,1];H^1_{\frac{2}{2-\alpha_0}})} \lesssim \|v\|_{X}^\alpha.
\end{equation}
Similarly,
\begin{equation}\label{E:d2A2npf2}
\begin{aligned}
	\| A_{2,n}(\cdot,v_1)-A_{2,n}(\cdot,v_2) \|_{L^{\frac{2}{\alpha_0-1}}((0,1];H^1_{\frac{2}{2-\alpha_0}})}
	&\lesssim
	\|F_n(v_1)-F_n(v_2)\|_{L^\infty((0,1];L^{\frac2{2-\alpha_0}})} \\
	&\lesssim n(\|v_1\|_X+\|v_2\|_X)^{\alpha-1} \|v_1-v_2\|_X
\end{aligned}
\end{equation}
for $v_1,v_2\in Z_M$.
For $\alpha_0 \in [1,\frac{7}{5}]$, we choose $a=3-2\alpha_0$.
By fractional Hardy and Gagliardo-Nirenberg inequality, one obtains
\[
	\| |x|^{-\frac{3-2\alpha_0}{\alpha_0+1}} v(t) \|_{L^{\frac{2 (\alpha_0+1)}{2-\alpha_0}}} \lesssim
	\| p^{\frac{3-2\alpha_0}{\alpha_0+1}} v(t) \|_{L^{\frac{2(\alpha_0+1) }{2-\alpha_0}}}
	\lesssim \|v(t)\|_{L^2}^{\frac{\alpha-2}{\alpha}}
\|p v(t)\|_{L^2}^{\frac{2}{\alpha}}.
\]
Hence,
\begin{equation}\label{E:d2A2npf3}
	\|A_{2,n}\|_{L^{\frac2{\alpha_0-1}} ((0,1];H^1_{\frac2{2-\alpha_0}})}
\lesssim  \|v\|_{L^\infty ((0,1];H^1)}^\alpha.
\end{equation}
Similarly,
\begin{equation}\label{E:d2A2npf4}
\begin{aligned}
	&{}\| A_{2,n}(\cdot,v_1)-A_{2,n}(\cdot,v_2) \|_{L^{\frac{2}{\alpha_0-1}}((0,1];H^1_{\frac{2}{2-\alpha_0}})} \\
	\lesssim{}&
	\||x|^{-3+2\alpha_0}(F_n(v_1)-F_n(v_2))\|_{L^\infty((0,1];L^{\frac2{2-\alpha_0}})} \\
	\lesssim{}& n(\|v_1\|_X+\|v_2\|_X)^{\alpha-1} \|v_1-v_2\|_X
\end{aligned}
\end{equation}
for $v_1,v_2 \in Z_M$.
Next, we estimate $A_{2,n}$ in $L^\infty H^1$.
One sees from
\begin{equation}\label{E:d2A3pf1}
	\left\| \J{\sqrt{t}p} U(t) R_n(t) ^{\frac12}U(-t) \cdot U(t) R_n(t) ^{\frac{\alpha_0-1}2}U(-t) \cdot U(t) R_n(t) ^{\frac{2-\alpha_0}2}V_n(t) \left|\frac{ x}{\sqrt{t}}\right|^{2-\alpha_0} \right\|_{\mathcal{B}(L^r)} \lesssim n^{\frac{\alpha_0-2}2}
\end{equation}
for $r=2$
that
\begin{align*}
	\| A_{2,n}(t) \|_{H^1}
	\lesssim n^{\frac{\alpha_0-2}2}
 t^{\alpha_0-1}\left\|  \J{p} \J{\sqrt{t} p}^{-1}  \left|\frac{x}{\sqrt{t}}\right|^{\alpha_0-2} F_n(v(t)) \right\| _{L^{2}}
 \lesssim n^{\frac{\alpha_0-2}2}
t^{\frac{\alpha_0-1}2}\| |x|^{-\frac{2-\alpha_0}{\alpha_0+1}} v(t) \|^{\alpha}_{L^{2(\alpha_0+1)}}.
\end{align*}
By the fractional Hardy inequality and Sobolev inequality, one obtains
\[
	\| |x|^{-\frac{2-\alpha_0}{\alpha_0+1}} v(t) \|_{L^{2(\alpha_0+1)}} \lesssim
	\| p^{\frac{2-\alpha_0}{\alpha_0+1}} v(t) \|_{L^{2(\alpha_0+1)}} \lesssim
 \|v(t)\|_{L^2}^{\frac{\alpha-2}{\alpha}}
\|p v(t)\|_{L^2}^{\frac{2}{\alpha}}.
\]
Hence, we find that
\begin{equation}\label{E:d2A2npf5}
	\|A_{2,n}(\cdot,v)\|_{L^\infty((0,1];H^1)} \lesssim n^{\frac{\alpha_0-2}2} \|v\|_{L^\infty((0,1];H^1)}^\alpha.
\end{equation}
A similar argument shows
\begin{equation}\label{E:d2A2npf6}
\begin{aligned}
	\| A_{2,n}(\cdot,v_1)-A_{2,n}(\cdot,v_2) \|_{L^\infty((0,1];H^1)}
	&\lesssim
	n^{\frac{\alpha_0-2}2}\||x|^{-2+\alpha_0}(F_n(v_1)-F_n(v_2))\|_{L^\infty((0,1];L^{2})} \\
	&\lesssim n^{\frac{\alpha_0}2}(\|v_1\|_X+\|v_2\|_X)^{\alpha-1} \|v_1-v_2\|_X
\end{aligned}
\end{equation}
for $v_1,v_2 \in Z_M$.
Thus,
we see from \eqref{E:d2A2npf1}, \eqref{E:d2A2npf3}, and \eqref{E:d2A2npf5} that
\begin{equation}\label{E:d2A2n1}
	\| A_{2,n} (\cdot,v)\|_X
	\lesssim  M^\alpha.
\end{equation}
Further, \eqref{E:d2A2npf2}, \eqref{E:d2A2npf4}, and \eqref{E:d2A2npf6} give us
\begin{equation}\label{E:d2A2n2}
	\| A_{2,n} (\cdot,v_1) - A_{2,n} (\cdot,v_2)\|_X
	\lesssim n M^{\alpha-1} \|v_1-v_2\|_X
\end{equation}
for $v_1,v_2 \in Z_M$.

We turn to the estimate of $A_{3,n}$ given in \eqref{E:A3def}. The Strichartz estimate implies that
\[
	\|A_{3,n}(\cdot,v)\|_{X} \lesssim
 \left\| U(s)
	\left( \frac{d}{ds}  R_n(s)\right) s^{\alpha_0-1} V_n(s) F_n(v(s))
\right\|_{L^{q_0'}((0,1];H^1_{r_0'})}.
\]
By means of Lemma \ref{L:Rest2} and Lemma \ref{L-K5/4-1},
one sees that
\[
	\left\| \J{\sqrt{s}p} U(s) R_n(s) ^{\frac12} U(-s)\right\|_{\mathcal{B}(L^{\frac2{\alpha_0}})} \lesssim 1,
\]
\[
	\left\|  U(s)R_n(s) ^{-\frac12}\left( \frac{d}{ds} R_n(s)\right) R_n(s)^{\frac{\alpha_0-3}4}U(-s)  \right\|_{\mathcal{B}(L^{\frac2{\alpha_0}})} \lesssim  n^{\frac{9-\alpha_0}8}
 s^{-1},
\]
and
\[
	\left\| U(s) R_n(s)^{\frac{3-\alpha_0}4} V_n(s) \left|\frac{ x}{\sqrt{s}}\right|^{\frac{3-\alpha_0}2} \right\|_{\mathcal{B}(L^{\frac2{\alpha_0}})} \lesssim n^{\frac{\alpha_0-3}4}.
\]
Hence, estimating as in the proof of \eqref{E:d2A2npf3}, we obtain
\begin{align}
	\begin{aligned}
	\| A_{3,n}(\cdot,v)\|_{X} \lesssim{}& n^{\frac{\alpha_0+3}8}
	\left\| s^{\frac{3\alpha_0-7}4} |x|^{\frac{\alpha_0-3}2} F_n(v) \right\| _{L^{q_0'}((0,1];L^{r_0'})} \\
	\lesssim{}& n^{\frac{\alpha_0+3}8}
	\|v\|_{L^\infty ((t,1],H^{1})}^\alpha
	\le n^{\frac{\alpha_0+3}8}M^\alpha
	\end{aligned}\label{E:d2A3n1}
\end{align}
since $s^{\frac{3\alpha_0-7}4} \in L^{q_0'}((0,1])$ $\Leftrightarrow$ $\alpha_0>1$ and
\[
	\| |x|^{\frac{\alpha_0-3}{2(\alpha_0+1)}} v \|_{L^{\frac{2(\alpha_0+1)}{\alpha_0}}} \lesssim
	\| |p|^{\frac{3-\alpha_0}{2(\alpha_0+1)}} v \|_{L^{\frac{2(\alpha_0+1)}{\alpha_0}}}
	\lesssim \| v \|_{L^2}^{\frac{3(\alpha_0-1)}{2(\alpha_0+1)}} \| \nabla v \|_{L^2}^{\frac{5-\alpha_0}{2(\alpha_0+1)}}
\]
for $\alpha_0>1$.
Similarly, for $\alpha_0>1$,
\begin{equation}\label{E:d2A3n2}
\begin{aligned}
	\| A_{3,n}(\cdot,v_1) -A_{3,n}(\cdot,v_2)\|_{X}
	&\lesssim n^{\frac{\alpha_0+3}8}
 \left\| |x|^{\frac{\alpha_0-3}2} (F_n(v_1)-F_n(v_2)) \right\| _{L^{\infty}((0,1];L^{r_0'})} \\
	&\lesssim n^{\frac{\alpha_0+11}{8}}
M^{\alpha-1} \|v_1-v_2\|_X.
\end{aligned}
\end{equation}

We move on to the estimate of $A_{4,n}(t)$.
By Strichartz,
\begin{align} \label{K8/21-1}
	\begin{aligned}
	\|A_{4,n}(\cdot,v)\|_{X} &\lesssim \left\|
	R_n(s) s^{\alpha_0-1} V_n(s) \frac{\partial F_n}{\partial z} (v(s)) p^2 v(s)
	\right\|_{L^1 ((0,1];H^1)}\\
	&\quad +  \left\| R_n(s) s^{\alpha_0-1} V_n(s) \frac{\partial F_n}{\partial \overline{z}} (v(s)) \overline{p^2 v(s)}
	\right\|_{L^1 ((0,1];H^1)}.
	\end{aligned}
\end{align}
We only consider the first term
since the both terms are handled similarly.
One sees from Lemma \ref{L:Rest2} that
\begin{align} \label{K8/21-2}
	\left\| \J{\sqrt{s}p} U(s)R_n(s) ^{\frac12} U(-s)\cdot U(s)  R_n(s)^{\frac{1}2} V_n(s) \J{\sqrt{s} p}
	\right\|_{\mathcal{B}(L^{2})} \lesssim n^{\frac{1}2},
\end{align}
which yields
\[
	\left\| U(s)
	 R_n(s) s^{\alpha_0-1} V_n(s) \frac{\partial F_n}{\partial z} (v(s))
	p^2 v(s)
\right\|_{H^1}
\lesssim n^{\frac12} s^{\alpha_0-2} \left\| \J{p}^{-1} \frac{\partial F_n}{\partial z} (v(s)) p^2 v(s) \right\| _{L^2}.
\]
We claim that
\begin{align}\label{E:d2A4nclaim}
	\left\|   \J{p}^{-1} \frac{\partial F_n}{\partial z} (v(s)) p^2 v(s) \right\| _{L^2}
	\lesssim n^2 \|v\|_{L^2}^{\alpha-2} \norm{v}_{H^{1}_{r_{0}}}
\| v\|_{H^{1}_{\frac{4}{3-\alpha_0} }  } + n\|v\|^{\alpha-1}_{H^1}\| v\|_{H^{1}_{\frac{4}{3-\alpha_0} }  }.
\end{align}
Indeed, by duality,
\begin{align*}
		\left\|   \J{p}^{-1} \frac{\partial F_n}{\partial z} (v(s)) p^2 v(s) \right\| _{L^2}
	&=  \sup_{\|\varphi\|_{L^2}=1}
	\left| \left(  p v,  p \left( \overline{\frac{\partial F_n}{\partial z} (v)} \J{p}^{-1}\varphi \right) \right)\right| \\
& \lesssim \|v\|_{H^1_{\frac4{3-\alpha_0}}}
\sup_{\|\varphi\|_{L^2}=1}
\left\| p \left( \overline{\frac{\partial F_n}{\partial z} (v)} \J{p}^{-1}\varphi \right) \right\|_{L^{\frac4{\alpha_0+1}}}.
\end{align*}
By H\"older,
\begin{align*}
	\sup_{\|\varphi\|_{L^2}=1}
\left\|   \overline{\frac{\partial F_n}{\partial z} (v)} (p\J{p}^{-1}\varphi)  \right\|_{L^{\frac4{\alpha_0+1}}} &{}\lesssim n \left\| v \right\|_{L^{\frac{4\alpha_0}{\alpha_0-1}}}^{\alpha-1}  \lesssim n\|v\|_{H^1}^{\alpha-1}.
\end{align*}
By using an estimate for the product of functions in Triebel-Lizorkin space \cite[\S 4.4.4 Theorem 3]{RS}, one obtains
\[
	\left\|   \overline{\left(p\frac{\partial F_n}{\partial z} (v)\right)} \J{p}^{-1}\varphi  \right\|_{L^{\frac4{\alpha_0+1}}}
\lesssim_\ep \left\| p\frac{\partial F_n}{\partial z} (v) \right\|_{L^{\frac4{\alpha_0-1+\ep}}}
\|\J{p}^{-1} \varphi\|_{H^1}
\]
for any $\ep \in (0,2)$.
We choose $\ep = 3-\alpha_0\in (0,2)$ to get
\[
	\left\| p\frac{\partial F_n}{\partial z} (v) \right\|_{L^{\frac4{\alpha_0-1+\ep}}}
\lesssim n^2 \|v\|_{H^1_{r_{0}}} \|v\|_{L^{2}}^{\alpha-2} \lesssim
n^2 \|v\|_{H^1}^{\alpha-2} \|v\|_{H^1_{r_{0}}}.
\]
We obtain the claim \eqref{E:d2A4nclaim}.
By interpolation, one has
\[
	\|v\|_{H^1_{\frac4{3-\alpha_0}}} \lesssim \|v\|_{H^1}^\frac12
	\|v\|_{H^1_{r_0}}^\frac12.
\]
Hence,
\begin{align*}
	\|A_{4,n}(\cdot,v)\|_{X} &\lesssim n^{\frac52} \|s^{\alpha_0-2} \|_{L^{(2q_0/3)'}((0,1])} \|v\|_{L^\infty((0,1];H^1)}^{\alpha-\frac{3}{2}} \|v\|_{L^{q_0}((0,1];H^1_{r_0})}^{\frac{3}{2}} \\
	&\quad + n^{\frac32}  \|s^{\alpha_0-2} \|_{L^{(2q_0)'}((0,1])} \|v\|_{L^\infty((0,1];H^1)}^{\alpha-\frac12} \|v\|_{L^{q_0}((0,1];H^1_{r_0})}^{\frac12}.
\end{align*}
Since $s^{\alpha_0-2} \in L^{(2q_0/3)'}((0,1]) \cap L^{(2q_0)'}((0,1])$ for $\alpha_0>1$, one obtains
\begin{equation}\label{E:d2A4n1}
	\|A_{4,n}(\cdot,v)\|_{X} \lesssim n^{\frac52} M^{\alpha-1}.
\end{equation}
Mimicking this argument, one also has
\begin{align*}
	\|A_{4,n}(\cdot,v_1)-A_{4,n}(\cdot,v_2)\|_{X}
	\lesssim{}& n^\frac12
\left\| p\frac{\partial F_n}{\partial z} (v_1) \right\|_{L^{q_{0}}((0,1];L^2)} \|v_1-v_2\|_{L^{2q_{0}}((0,1];H^1_{\frac{4}{3-\alpha_0}})} \\
 	&+n^{\frac{1}{2}} \left\| \frac{\partial F_n}{\partial z}(v_1) \right\|_{L^{\infty}((0,1];L^{\frac4{\alpha_{0}-1}})} \|v_{1} - v_{2}\|_{L^{2q_0}((0,1];H^1_{\frac{4}{3-\alpha_0}})} \\
	&+n^\frac12\left\| \frac{\partial F_n}{\partial z} (v_1) -\frac{\partial F_n}{\partial z} (v_2) \right\|_{L^{\infty}((0,1];L^{\frac4{\alpha_0-1}})} \|v_2\|_{L^{2q_{0}}((0,1];H^1_{\frac{4}{3-\alpha_0}})}
\\
&+n^\frac12 \left\| p\left(\frac{\partial F_n}{\partial z} (v_1) -\frac{\partial F_n}{\partial z} (v_2) \right) \right\|_{L^{q_{0}}((0,1];L^2)} \|v_2\|_{L^{2q_{0}}((0,1];H^1_{\frac{4}{3-\alpha_0}})}
\\
\lesssim{}& n^{\frac{7}2}(\|v_1\|_X+\|v_2\|_X)^{\alpha-1} \|v_1-v_2\|_X
\end{align*}
for $\alpha_0>1$. Hence,
\begin{equation}\label{E:d2A4n2}
	\|A_{4,n}(\cdot,v_1)-A_{4,n}(\cdot,v_2)\|_{X} \lesssim n^{\frac{7}2}M^{\alpha-1} \|v_1-v_2\|_X
\end{equation}
for $v_1,v_2 \in Z_M$.

We lastly consider the estimate of $A_{5,n}$ given in \eqref{E:A5def}.
By Strichartz,
\begin{align*}
	\|A_{5,n}\|_{X} \lesssim {}&
 \left\| U(s)
	 R_n(s) s^{2\alpha_0-3} V_n(s) \frac{\partial F_n}{\partial z} (v(s)) \tilde{F}(s,v(s)) \right\|_{L^{q_0'} ((0,1];H^1_{r_0'})}\\
 &+\left\| U(s)
	 R_n(s) s^{2\alpha_0-3} V_n(s) \frac{\partial F_n}{\partial \overline{z}} (v(s)) \overline{\tilde{F}(s,v(s))} \right\|_{L^{q_0'} ((0,1];H^1_{r_0'})}.
\end{align*}
Let us treat only the first term.
We use the bound \eqref{E:d2A3pf1} with $r=2/\alpha_0$ to obtain
\begin{multline*}
	\left\| U(s)
	 R_n(s) s^{2\alpha_0-3} V_n(s)
\frac{\partial F_n}{\partial z} (v(s)) \tilde{F}(s,v(s)) \right\|_{H^1_{r_0'}}\\
	\lesssim
n^{\frac{\alpha_0-2}2}
 s^{\frac{3\alpha_0-5}2}
\left\| |x|^{\alpha_0-2} \frac{\partial F_n}{\partial z} (v(s)) \tilde{F}(s,v(s)) \right\|_{L^{\frac{2}{\alpha_0}}}.
\end{multline*}
By fractional Hardy inequality,
one finds that
\begin{align*}
	\left\| |x|^{\alpha_0-2} \frac{\partial F_n}{\partial z} (v(s))\tilde{F}(s,v) \right\|_{L^{\frac{2}{\alpha_0}}}
&\lesssim n\|\lambda_m\|_{\ell^1_m} \| |x|^{-\frac{2-\alpha_0}{2\alpha_0+1}} v\|_{L^{\frac{2(2\alpha_0+1)}{\alpha_0} }}^{2\alpha-1}\\
&\lesssim n\|\lambda_m\|_{\ell^1_m} \| |p|^{\frac{2-\alpha_0}{2\alpha_0+1}} v\|_{L^{\frac{2(2\alpha_0+1)}{\alpha_0} }}^{2\alpha-1}.
\end{align*}
By Gagliardo-Nirenberg inequality, we further obtain
\[
	\| |p|^{\frac{2-\alpha_0}{2\alpha_0+1}} v\|_{L^{\frac{2(2\alpha_0+1)}{\alpha_0} }} \lesssim \| v\|_{L^2}^{\frac{2(\alpha_0-1)}{2\alpha_0+1}} \| p v\|_{L^2}^{\frac{3}{2\alpha_0+1}}.
\]
Thus, we end up with the estimate
\begin{equation}\label{E:d2A5n1}
	\|A_{5,n}(\cdot,v)\|_{X} \lesssim
	n^{\frac{\alpha_0}2}\|\lambda_m\|_{\ell^1_m} \|v\|_{L^\infty ((0,1];H^1)}^{2\alpha-1} \le n^{\frac{\alpha_0}2}M^{2\alpha-1}.
\end{equation}
A similar treatment yields
\begin{align*}
	\|A_{5,n}(\cdot,v_1)-A_{5,n}(\cdot,v_2)\|_{X} &\lesssim  n^{\frac{\alpha_0}2}
\left\| |x|^{\alpha_0-2} \left(\frac{\partial F_n}{\partial z} (v_1) \tilde{F}(\cdot,v_1)-
\frac{\partial F_n}{\partial z} (v_2) \tilde{F}(\cdot,v_2)\right) \right\|_{L^\infty ((0,1];L^{\frac{2}{\alpha_0}})} \\
&\lesssim  n^{\frac{\alpha_0}2+1}
(\|v_1\|_X+\|v_2\|_X)^{2\alpha-2} \|v_1-v_2\|_X
\end{align*}
for $v_1,v_2\in Z_M$.
Hence, we obtain
\begin{equation}\label{E:d2A5n2}
	\|A_{5,n}(\cdot,v_1)-A_{5,n}(\cdot,v_2)\|_{X} \lesssim n^{\frac{\alpha_0}2+1} M^{2(\alpha-1)} \|v_1-v_2\|_X
\end{equation}
for $v_1,v_2\in Z_M$.

The estimate
\eqref{E:Inest1} for $d=2$
follows from
\eqref{E:d2A1n1}, \eqref{E:d2A2n1},
\eqref{E:d2A3n1},
\eqref{E:d2A4n1},
and \eqref{E:d2A5n1}.
Similarly, \eqref{E:d2A2n2},
\eqref{E:d2A3n2},
\eqref{E:d2A4n2},
and \eqref{E:d2A5n2} yield \eqref{E:Inest2}.

\subsubsection{Proof of Lemma \ref{L:Inest} in three dimensions}
We next consider the three dimensional case of the first assertion of Lemma \ref{L:Inest}.
Recall that
\[
	(q_0,r_0) = \left( 4,3 \right), \quad \alpha_0 = \frac32, \quad \alpha =2.
\]
Here, we remark that the Strauss exponent for $d=3$ equals to $2$, and hence we aim to establish GWP of \eqref{1} when $\alpha \leq 2$. On the other hand, to discuss the GWP for \eqref{E:v1} in the $H^1$ framework, we need to handle the term
\begin{align}\label{K1/14-1}
|p| R_n(s) V_n(s) (\partial_z F_n) p^2 v(s),
\end{align}
which arises in the estimation to $A_{4,n} (t)$.
Since $|p|R_n(s)|p|$ is bounded, we extract $|p|$ from $p^2=|p||p|$ and apply it to $|p|R_n(s)$ to obtain the boundedness for \eqref{K1/14-1}.
The commutator $[\partial_z F_n, |p|]$ appears in this process, and hence $\partial _z F_n$ must be differentiable at least once to justify the commutator.
This imposes the restriction $\alpha -1 \geq 1$.
Consequently, our method is applicable only for the case $\alpha =2$.

Fix $n\ge2$ and pick $v,v_1,v_2\in Z_M$. The calculation is similar to the case where $d=2$ and hence, we give a sketch of the estimates related to $A_{j,n} (t)$ for $j=2$, $3$, $4$, $5$.

Let us first estimate $A_{2,n}$ given in \eqref{E:A2def}.
Arguing as in the case $d=2$, we see from the fractional Hardy equality and Gagliardo-Nirenberg equality that
\begin{align*}
	\| A_{2,n}(t,v) \|_{H^1_{} }
	\lesssim n^{- \frac14}t^{\frac34} \left\|  \J{p} \J{\sqrt{t} p}^{-1} \right\|_{\mathcal{B}(L^{2})}  \left\| |x|^{-\frac12}   F_n(v(t)) \right\| _{L^{2}} &  \lesssim n^{- \frac14}t^{\frac14}  \| |x|^{- \frac14} v(t) \|_{L^4} ^2
	\\
	& \lesssim n^{- \frac14}t^{\frac14}  \| v(t) \|_{H^1} ^2,
\end{align*}
and we see from the Gagliardo-Nirenberg equality that
\begin{align*}
	\| A_{2,n}(t,v) \|_{H^1_{3} }
	\lesssim t^{\frac12} \left\|  \J{p} \J{\sqrt{t} p}^{-1} \right\|_{\mathcal{B}(L^{3})}  \left\|  F_n(v(t)) \right\| _{L^{3}}
	\lesssim   \|  v(t) \|^{2}_{L^{6}} \lesssim \| v(t) \|_{H^1} ^2.
\end{align*}
Hence, for $(\tau, \theta ) = ( \infty, 2)$ or $(4, 3)$,
\begin{align} \label{dim3:1}
	\|A_{2,n}\|_{L^{\tau} ((0,1];H^1_{\theta})}
	\lesssim \norm{t^{\frac{1}{4} - \frac1{\tau} } }_{L^{\tau}((0,1])} \|v\|_{L^\infty ((0,1];H^1)}^2.
\end{align}
Thus, one has
\begin{align*}
	\| A_{2,n} (\cdot,v)\|_X
	\lesssim  M^{2}.
\end{align*}
Also, As for the estimate of $A_{3,n}$ given in \eqref{E:A3def}, we have
\begin{align*}
	s^{- \frac12} \left\| |\sqrt{s} p| U(s) \left( \frac{d}{ds}  R_n(s)\right) s^{\alpha_0-1} V_n(s) F_n(v(s)) \right\|_{L^{r_0'}}
	&\lesssim n^{\frac12} s^{-1+\frac{1}{2}} \norm{\abs{{x}{}}^{-{1}}F_{n}(v(s)) }_{L^{3/2}} \\
	&\lesssim n^{\frac12} s^{- \frac12} \|v(s)\|_{H^1}^2.
\end{align*}
By the Strichartz estimate, this implies
\begin{align} \label{dim3:2}
	\| A_{3,n}(\cdot,v)\|_{X} \lesssim n^{\frac12} \norm{s^{- \frac12}}_{L^{4/3}_{s}((0,1])} \|v\|_{L^\infty ((0,1];H^1)}^2 \le n^{\frac12} M^2.
\end{align}

Let us move on to the estimate of $A_{4,n}(t)$.
We handle the first term of right hand side of \eqref{K8/21-1}.
One sees from \eqref{K8/21-2} that
\[
	\left\| U(s)
	 R_n(s) s^{\frac12} V_n(s) \frac{\partial F_n}{\partial z} (v(s))
	p^2 v(s)
\right\|_{H^1}
\lesssim n^{\frac12} s^{-\frac12} \left\| \J{p}^{-1} \frac{\partial F_n}{\partial z} (v(s)) p^2 v(s) \right\| _{L^2}.
\]
We claim that
\begin{align*}
	\left\|   \J{p}^{-1} \frac{\partial F_n}{\partial z} (v(s)) p^2 v(s) \right\| _{L^2}
	\lesssim (n+n^2) \| v\|_{H^{1}} \| v \|_{H_3^1}  .
\end{align*}
Indeed, by duality,
\begin{align*}
		\left\| \J{p}^{-1} \frac{\partial F_n}{\partial z} (v(s)) p^2 v(s) \right\| _{L^2}
	&=  \sup_{\|\varphi\|_{L^2}=1}
	\left| \left(  p v,  p \left( \overline{\frac{\partial F_n}{\partial z} (v)} \J{p}^{-1}\varphi \right) \right)\right| \\
& \lesssim \|v\|_{H^1_{3}}
\sup_{\|\varphi\|_{L^2}=1}
\left\| p \left( \overline{\frac{\partial F_n}{\partial z} (v)} \J{p}^{-1}\varphi \right) \right\|_{L^{\frac32}}.
\end{align*}
By H\"older,
\begin{align*}
	\sup_{\|\varphi\|_{L^2}=1}
\left\|   \overline{\frac{\partial F_n}{\partial z} (v)} (p\J{p}^{-1}\varphi)  \right\|_{L^{\frac32}} &{}\lesssim n \left\| v \right\|_{L^{6}} \lesssim n\|v\|_{H^1}.
\end{align*}
Moreover by the H\"{o}lder,
\[
	\left\|   \overline{\left(p\frac{\partial F_n}{\partial z} (v)\right)} \J{p}^{-1}\varphi  \right\|_{L^{\frac32}}
\lesssim \left\| p\frac{\partial F_n}{\partial z} (v) \right\|_{L^{2}}
\|\J{p}^{-1} \varphi\|_{L^6} \lesssim n^2 \| v \|_{H^1}.
\]
Hence,
\begin{align} \label{dim3:3}
	\|A_{4,n}(\cdot,v)\|_{X} &\lesssim n^{\frac{1}{2}} (n+n^2) \left( \| v \|_{L^{\infty}((0,1];H^1) }  \| v \|_{L^{4}((0,1];H^1_3)}  \right) \lesssim n^{\frac{5}{2}} M^2,
\end{align}
where we have used $ \left\| s^{- \frac12} \| v(s) \|_{H_3^1} \right\|_{L^1_s((0,1])} \lesssim \norm{s^{-\frac{1}{2}}}_{L^{4/3}_{s}} \| v \| _{L^{4}((0,1];H^1_3)} \lesssim \| v \| _{L^{4}((0,1];H^1_3)}$.

We lastly consider the estimate of $A_{5,n}$ given in \eqref{E:A5def}.
Using the Strichartz estimate,
\begin{align*}
	\|A_{5,n}(\cdot, v)\|_{X} &\lesssim
 \left\| U(s)
	 R_n(s)  V_n(s) \frac{\partial F_n}{\partial z} (v(s)) \tilde{F}(s,v(s)) \right\|_{L^{1} ((0,1];H^1_{})} \\
 	&\quad + \left\| U(s)
	 R_n(s)  V_n(s) \frac{\partial F_n}{\partial \overline{z}} (v(s)) \overline{\tilde{F}(s,v(s))} \right\|_{L^{1} ((0,1];H^1_{})}
	 \\
	 &\lesssim n \|\lambda_m\|_{\ell^1_m} \left\| s^{- \frac12}  \left\| v(s) \right\|_{L^6} ^3 \right\|_{L^1_s((0,1])} \\
	 &\lesssim \norm{s^{-\frac{1}{2}}}_{L^{1}_{s}} n \| v \|_{L^{\infty} ((0,1];H^1)}^3
	 \lesssim n \| v \|_{L^{\infty} ((0,1];H^1)}^3,
\end{align*}
and thus we end up with the estimate
\begin{align} \label{dim3:4}
	\|A_{5,n}(\cdot,v)\|_{X} \lesssim
	n  \|v\|_{L^\infty ((0,1];H^1)}^{3} \le n^{}M^{3}
\end{align}
for $v_1,v_2\in Z_M$.

\subsection{Proof of Proposition \ref{P:key}}\label{ss:keyprop}

We are now in the position to prove our key proposition.
\begin{proof}[Proof of  Proposition \ref{P:key}]
Let us consider the first assertion.
Pick $v\in Z_M$.
Plugging \eqref{E:XtoY},\eqref{E:Pkeypf1}, \eqref{E:Pkeypf2}, \eqref{E:Pkeypf3}, and \eqref{E:Inest1},
one sees that there exists $C>0$ independent of $T$ such that
\[
	\|\Phi(v)\|_{X} + \|\Phi(v)\|_{Y}
\le C(\ep_0+\ep_0^\alpha) + C\| \J{n}^\frac52 \lambda_n\|_{\ell^1}M^\alpha(1+M^{\alpha-1}).
\]
We fix $M>0$ is so small that
\[
	C\| \J{n}^\frac52 \lambda_n\|_{\ell^1}M^\alpha(1+M^{\alpha-1})
	\le \frac12 M
\]
and take $\ep_0$ so small that
\[
	C(\ep_0+\ep_0^\alpha) \le \frac12 M.
\]
Under this choice, we have $\|\Phi(v)\|_X \le M$, i.e.,
$\Phi(v) \in Z_M$.
Hence, $\Phi$ is a map from $Z_M$ to itself.

To see $\Phi$ is a contraction,
we combine \eqref{E:Pkeypf2.5}, \eqref{E:Pkeypf4}, and \eqref{E:Inest2} to obtain
\[
	\|\Phi(v_1)-\Phi(v_2)\|_X \le C\| \J{n}^\frac{7}{2} \lambda_n\|_{\ell^1} M^{\alpha-1} \J{M}^{\alpha-1} \|v_1-v_2\|_X
\]
for $v_1,v_2\in Z_M$. Letting $M$ (and $\ep_0$) even smaller if necessary, one sees that
\[
	\|\Phi(v_1)-\Phi(v_2)\|_X \le \frac12 \|v_1-v_2\|_X,
\]
which completes the proof.

We turn to the proof of the second assertion.
Let $v\in Z_M$ be a fixed point of \eqref{E:IE}.
By \eqref{E:Phierror} and the fact that $v=\Phi(v)$, one has
\[
	\CAL{E}(v) = \sum_{n\ge2} \lambda_n t^{\alpha_0-1} U(t) R_n(t) V_n(t) \left( \frac{\partial F_n}{\partial z} \CAL{E}(v) - \frac{\partial F_n}{\partial \overline{z}} \overline{\CAL{E}(v)} \right).
\]
Mimicking the estimate in the proof of Lemma \ref{L:Phierror},
one sees that
\[
	\| \CAL{E}(v) \|_{L^{q_1}((0,1];H^{-1}_{r_1})} \lesssim M^{\alpha-1} \| \CAL{E}(v) \|_{L^{q_1}((0,1];H^{-1}_{r_1})}.
\]
If $M$ is small then this inequality implies $\CAL{E}(v)=0$.
This implies that $v$ is a solution to \eqref{E:Iv1}.
\end{proof}

\subsubsection{Existence and uniqueness for a solution in Definition \ref{D:sol}. } \label{3.3.1}

Given suitable $u_0$, we obtain a solution $v(t)$
to \eqref{E:v1} on $(0,1]$, in the sense of Definition \ref{D:v1sol}, via Proposition \ref{P:key}.
Then, by means of Lemma \ref{L7},
we obtain a solution $u(t)$ to \eqref{1} on $[0,\infty)$ in the sense of Definition \ref{D:sol}.
Further, the solution scatters in $H^{0,1}$ in the sense of Theorem \ref{T1}.
Let us next prove the uniqueness of the solution to \eqref{1}.
To clarify the argument, let us
denote by
$u_*(t)$
 the solution to \eqref{1} obtained by the above procedure.
Similarly, let $v_*(t)$ be the above solution to \eqref{E:v1} on $(0,1]$.

Suppose $u(t)$ be a solution to \eqref{1} in the sense of Definition \ref{D:sol} for the same initial data $u_0$.
The goal is to prove that $u$ coincides with $u_*$.
First, we let
\begin{align*}
v(t,x) = t^{- \frac{d}2} e^{i \frac{x^2}{2t}} \overline{u\left( \frac{1-t}t, \frac{x}t \right) }.
\end{align*}
Then, one sees from the standard argument on the pseudo-conformal transformation
that $v(t)$ is a solution to \eqref{E:v1} in the sense of Definition \ref{D:v1sol}.
Hence, it suffices to show that $v=v_*$.
To this end, we shall prove that $v$ is a fixed point of $\Phi$
and belongs to $Z_M$, where $M=M(u_0)$ is the same number as in Proposition \ref{P:key}.
If we had these two properties, one  concludes $v=v_*$ from the uniqueness of the fixed point of $\Phi$ in $Z_M$.

Pick $T \in [0,1)$.
Let $Z_{M,T}$ be the same complete metric space as $Z_M$ but the time interval is replace by $(T,1]$.
Note that $Z_{M}=Z_{M,0}$ and $\Phi$ is a contraction map from $(Z_{M,T},d_X)$ to itself for any $T \in [0,1)$.
Since $v \in C((0,1];H^1) \cap L^{q_0}_{t,\mathrm{loc}}((0,1]; L^{r_0}_x)$,
there exists $T \in [0,1)$ such that $ v \in Z_{M,T}$, i.e.,
\[
	\| v \|_{L^\infty ((T,1];H^{1})} + \|v\|_{L^{q_0} ((T,1];H^{1}_{r_0})} + \|\partial_t v\|_{L^{q_1} ((T,1];H^{-1}_{r_1})} \le M.
\]
Let us prove that $v$ is a fixed point of $\Phi$.
Since $\Phi$ is a map from $Z_{M,T}$ to itself, one has $\Phi(v) \in Z_{M,T}$.
In particular,
since $\mathcal{E}(v) = 0$,
\begin{align*}
	\sum_{j=1}^5A_{j,n}(t)={}&-U(t)\int_{1}^{t}  \left( \frac{d}{ds}  R_n(s)\right) s^{\alpha_0-1} V_n(s) F_n(v(s)) ds\\
&\quad - U(t) \int_{1}^{t} R_n(s)s^{\alpha_0-1} V_n(s) \left( \frac{d}{ds} F_n(v(s)) \right) ds \\
&\quad + t^{\alpha_0-1} U(t) R_n(t) V_n(t) F_n(v(t)) - U(t) R_n(1)V_n(1)F_n(v(1))
\end{align*}
holds and further all terms in the right hand side
belongs to $C((T,1];H^1)$.

By an integration by parts, one sees that
\begin{align} \label{K7/16-1}
	I_n (t)  = \sum_{j=1}^5  A_{j,n}(t)
\end{align}
for all $n\ge2$.
Indeed, since integrals in $I_n$ and $A_{j,n}$ $(1\le j \le 5)$
make sense as a Bochner integral on the interval $[t,1]$
in $L^2$, they are differentiable
in the strong $L^2$-sense
almost every $t\in [T,1]$ (see e.g., \cite[Lemma 2.5.8]{YNML16}) and satisfy the identity
\begin{align*}
	\frac{d}{dt}\left( U(-t) \left( \sum_{j=1}^5 A_{j,n} (t) -I_n(t) \right)  \right) &=  -t^{\alpha_0 -2}  U(-t) e^{-i \frac{(n-1)x^2 }{2t}} F_n(v(t))
\\ & \quad  +   \frac{d}{dt} \left(   t^{\alpha_0 -1} R_n(t) V_n(t) F_n(v(t)) \right)
\\ & \qquad  - t^{\alpha_0 -1}  \left(  \frac{d}{dt} R_n(t) \right) V_n(t) F_n(v(t))
\\ &  \quad \qquad - t^{\alpha_0 -1} R_n(t) V_n(t) \frac{d}{dt} \left( F_n(v(t)) \right)
\\ &= 0.
\end{align*}
Together with $I_n(1)=0$ and $A_{j,n}(1)=0$ for $1\le j \le 5$, one has
\[
	U(-t) \left( \sum_{j=1}^5 A_{j,n} (t) -I_n(t) \right)
	=0+ \int_t^1
	\frac{d}{dt}\left( U(-t) \left( \sum_{j=1}^5 A_{j,n} (t) -I_n(t) \right)  \right)(s) ds=0
\]
for $t\in [T,1]$, from which \eqref{K7/16-1} follows.
Since the right hand side of \eqref{K7/16-1} belongs to $C([T,1];H^1)$ for all $n\ge 1$ thanks to \eqref{E:Pkeypf2} and Lemma \ref{L:Inest}, so does $I_n$.
Together with the weighted summability assumption on $\{ \lambda_n \}$, one sees that
\[
	-i \int_1^t s^{\alpha_0-2} U(t-s) \tilde{F}(s,v(s))ds = -i \sum_{n =1}^\infty  {\overline{\lambda_n}} I_n.
\]
In particular, $v$ satisfies \eqref{E:Iv1}.
Thus, one sees that
\begin{align*}
	\Phi(v) ={}& U(t-1) v(1) - i {\overline{\lambda_1}} I_1(t,v) -i \sum_{n =2}^\infty  {\overline{\lambda_n}} \sum_{j=1}^5A_{j,n} \\
	={}& U(t-1) v(1) -i \sum_{n =1}^\infty  {\overline{\lambda_n}} I_n \\
	={}& v(t)
\end{align*}
as desired.

Finally, we claim that
\[
	\inf \{ T \in [0,1)\ ;\ v =v_* \text{ on } (T,1]\} = 0.
\]
Note that we have shown in the previous paragraph that the left hand side is less than one.
To prove the claim by contradiction,
we suppose that the left hand side equals to a positive number $T_0 \in (0,1)$.
Then,
since $v \in Z_{M,0}$, and $v=v_*$ on $(T_0,1]$, we see that
\[
	\| v \|_{L^\infty ((T_0,1],H^{1})} + \|v\|_{L^{q_0} ((T_0,1],H^{1}_{r_0})} + \|\partial_t v\|_{L^{q_1} ((T_0,1],H^{-1}_{r_1})} <M
\]
and hence there exists $\tilde{T}_0\in [0,T_0)$ such that
$v\in Z_{M,\tilde{T}_0}$.
Arguing as in the previous paragraph, one deduces that $v$ is a fixed point of $\Phi$.
Then, one sees that $v=v_*$ on $(\tilde{T}_0,1]$.
However, this contracts with the definition of $T_0$.
Thus, we see that the claim is true.

The claim implies that $v \in Z_M$ and $v=v_*$ holds.
This is nothing but the desired uniqueness property.

\subsection{Proof of Scattering}\label{ss:scattering}
\subsubsection{Scattering for one time direction}
Here, we prove the  scattering-part of Theorem \ref{T1}.
By means of Lemma \ref{L7}, it suffices to prove that $v(t)$, given in \eqref{E:vdef}, has a limit as $t\to0+$.
We will prove that $\{U(-t)v(t)\}_{0<t\le1}$ has the Cauchy property as $t\to0+$.
Since $U(t)$ is unitary on $H^1$, it suffices to show that
\[
	\| v-U(\cdot -t_0)v(t_0)\|_{L^\infty((0,t_0];H^1)} \to 0
\]
as $t_0 \to 0+$.

Let $v\in Z_M$ be a solution to \eqref{E:IE}.
One has
\begin{align*}
	v(t) - U(t-t_0)v(t_0)
	&= -i{\overline{\lambda_1}} (I_1(t)-U(t-t_0)I_1(t_0)) \\
&\quad - i \sum_{n=2}^\infty {\overline{\lambda_n}} \sum_{k=2}^5 (A_{k,n}(t)-U(t-t_0)A_{k,n}(t_0)).
\end{align*}
Let us estimate each terms in the right hand side.

We begin with the term
\[
	I_1(t)-U(t-t_0)I_1(t_0)
	= \int_{t_0}^t s^{\alpha_0-2} U(t-s)  F_{1}( v(s) )
 ds .
\]
Mimicking the proof of \eqref{E:Pkeypf2} and \eqref{E:Pkeypf3}, one finds that
\[
	\|I_1-U(\cdot-t_0)I_1(t_0)\|_{L^\infty((0,t_0];H^1)}
	\lesssim
	\begin{cases}
	\|s^{\alpha_0-2}\|_{L^1_s((0,t_0])} M^\alpha
&d=1, \\
	\|s^{\alpha_0-2}\|_{L^{\frac{2}{3-\alpha_{0}}}_s((0,t_0])} M^\alpha
& d=2, 3.
	\end{cases}
\]
In the both cases, the right hand side tends to zero as $t_0\to0+$.

We next consider the estimate of $A_{2,n}(t)-U(t-t_0)A_{2,n}(t_0)$.
We estimate these two terms separately.
Mimicking the proof of \eqref{E:d1A2n1}, \eqref{E:d2A2npf5} and \eqref{dim3:1}, one has
\[
	\|A_{2,n}(t)\|_{L^\infty ((0,t_0];H^1)} \lesssim
\begin{cases}
 t_0^{\alpha_0-\frac54} M^\alpha & d=1 \\
t_0^{\frac{\alpha_0-1}2} M^\alpha & d=2, \\
t_0^{\frac14} M^\alpha & d=3.
\end{cases}
\]
The right hand side tends to zero as $t_0\to0+$.
Further, one has
\[
	\|U(t-t_0)A_{2,n}(t_0)\|_{L^\infty ((0,t_0];H^1)}
=\|A_{2,n}(t_0)\|_{H^1} \le
\|A_{2,n}(t)\|_{L^\infty ((0,t_0];H^1)}.
\]
This term is estimated in the exactly same way.

Let us move on to the estimate of
\[
	A_{3,n}(t)-U(t-t_0)A_{3,n}(t_0) = - U(t) \int_{t_0}^{t} s^{\alpha_0-1} \left( \frac{d}{ds}  R_n(s)\right)  V_n(s) F_n(v(s)) ds.
\]
Hence, estimating as in the proof of \eqref{E:d1A3n1}, \eqref{E:d2A3n1} and \eqref{dim3:2}, one obtains
\[
	\| A_{3,n}-U(\cdot-t_0)A_{3,n}(t_0) \|_{L^\infty ((0,t_0];H^1)} \lesssim
	\begin{cases}
	n \| s^{\alpha_0- \frac{9}{4}}\|_{L^{1}((0,t_0])}M^\alpha & d=1,\\
	n \| s^{\frac{3\alpha_{0}-7}{4}}\|_{L^{q_0'}((0,t_0])}M^\alpha & d=2, \\
	{n}^{\frac12} \| s^{-\frac12}\|_{L^{\frac43}((0,t_0])}M^\alpha & d=3.
	\end{cases}
\]
The right hand side tends to zero as $t_0\to0+$.

The last two term are handled similarly. Indeed, mimicking the proof of \eqref{E:d1A4n1}, \eqref{E:d2A4n1} and \eqref{dim3:3}, one obtains
\[
	\| A_{4,n}-U(\cdot-t_0)A_{4,n}(t_0) \|_{L^\infty ((0,t_0];H^1)} \lesssim
	\begin{cases}
	n^{\frac52} \|s^{\alpha_0-2}\|_{L^{\frac{12}{11}}((0,s_0])}M^\alpha &d=1,\\
	n^{\frac52} \|s^{\alpha_0-2}\|_{L^{(\frac{4}{3(\alpha_0-1)})'} ((0,t_0]) \cap L^{(\frac{4}{\alpha_0-1})'} ((0,t_0])}M^\alpha & d=2, \\
	n^{\frac52} \|s^{-\frac{1}{2}}\|_{L^{4/3} ((0,t_0])} M^\alpha & d=3,
	\end{cases}
\]
of which right hand sides tend to zero as $t_0\to0+$.
Further, estimates similar to \eqref{E:d1A5n1} and \eqref{E:d2A5n1} yield
\[
	\| A_{5,n}-U(\cdot-t_0)A_{5,n}(t_0) \|_{L^\infty ((0,t_0];H^1)} \lesssim
	\begin{cases}
	n \|s^{2\alpha_0-\frac{13}{4}}\|_{L^1((0,s_0])}M^{2\alpha-1} & d=1,\\
	n^{\frac{\alpha_0}2} \|s^{\frac{3\alpha_0-5}2}\|_{L^{(\frac{2}{\alpha_0-1})'} ((0,t_0])} M^{2\alpha-1} & d=2, \\
	n \|s^{-\frac{1}{2}}\|_{L^{1} ((0,t_0])} M^{2\alpha-1} & d=3.
	\end{cases}
\]
The right hand side tends to zero as $t_0\to0+$.

\subsubsection{Small data scattering in $H^{1,1}$ for both time directions}

In this subsection, we prove Corollary \ref{C:both}.

\begin{proof}[Proof of Corollary \ref{C:both}]
Let $\widetilde{\ep}_0 = (3 \cdot 2^{\frac{1}{\alpha-1} - \frac{d}{4}})^{-1} \ep_0$, where $\ep_0$ is the constant given in Theorem \ref{T1}.
Then, if $u_0 \in H^{1,1}$ satisfies $\|u_0\|_{H^{1,1}}\le \widetilde{\ep}_0$ then it satisfies the assumption \eqref{E:smallcond} for $t_0=\pm1,\pm2$.
Indeed,
\begin{align*}
	|t_0|^{\frac1{\alpha-1} - \frac{d}4}
	(\| u_0 \|_{L^2 }+|t_0|^{-\frac{1}2}\|J(t_0)u_0\|_{L^2})
	&\leq 2^{\frac{1}{\alpha-1} - \frac{d}{4}} \left( \norm{u_{0}}_{L^{2}} + \norm{xu_{0}}_{L^{2}} + \norm{pu_{0}}_{L^{2}} \right) \\
	&\leq 3 \cdot 2^{\frac{1}{\alpha-1} - \frac{d}{4}} \norm{u_{0}}_{H^{1, 1}} \leq \varepsilon_{0}
\end{align*}
if $|t_0| \le 2$.
Let us consider the scattering for the positive time direction.
By Theorem \ref{T1} with $t_0=1$ and $t_0=2$,
the solution exists on $[0,\infty)$ and satisfies $u(t) \in X^1(t+1) \cap X^1(t+2) = H^{1,1}$ for all $t\ge0$.
Further,
there exist $u_{+,1},u_{+,2} \in H^{0,1}$ such that
\[
	\| U(-t-j) u(t)-u_{+,j} \|_{H^{0,1}} = \| U(-t) u(t)-U(j)u_{+,j} \|_{X^1(j)}
	 \to 0
\]
as $t\to \infty$.
This shows $U(1)u_{+,1} = U(2) u_{+,2}$ in $L^2$.
We set $u_+:=U(1)u_{+,1}$.
Then, $u_+ \in X^1(1) \cap X^1(2) = H^{1,1}$.
Further, since $p=J(1)-J(2)$ and $x= 2J(1)-J(2)$
\begin{align*}
	\| U(-t)u(t) - u_+ \|_{H^{1,1}}
	\le{}& \| U(-t)u(t) - u_+ \|_{L^2} +
	\| p(U(-t)u(t) - u_+) \|_{L^2} \\&+ \| x(U(-t)u(t) - u_+) \|_{L^2} \\
	\le {}& \| U(-t)u(t) - u_+ \|_{L^2} +
	3\| J(1)(U(-t)u(t) - u_+) \|_{L^2} \\&+ 2\| J(2)(U(-t)u(t) - u_+) \|_{L^2} \\
	= {}& \| U(-t)u(t) - u_+ \|_{L^2} +
	3\| U(-t)u(t) - u_+ \|_{X^1(1)} \\&+ 2\| U(-t)u(t) - u_+ \|_{X^1(2)} \to 0
\end{align*}
as $t\to\infty$. Hence, scattering holds in $H^{1,1}$. The statement for the negative time direction is shown in the same way.
\end{proof}

\section*{Acknowledgments}
M.K. was supported by JSPS KAKENHI Grant Numbers 24K06796 and 20K14328.
S.M. was supported by JSPS KAKENHI Grant Numbers 24K00529, 23K20803 and 23K20805.
H.M. was supported by JSPS KAKENHI Grant Numbers 22K13941.
The authors would like to express their sincere gratitude to the anonymous referee for carefully reviewing their manuscript and providing valuable suggestions.

 ~~ \\   {\bf \Large Declarations} \\ ~~ \\ 
{\bf Conflict of interest } \\ ~~ \\  The authors declare that there is no conflict of interest.
\\ ~~ \\ 
{\bf Data Availability Statement } \\ ~~ \\ 
No data are analyzed in this study


\end{document}